\documentclass[10pt]{article}
\usepackage{amsfonts,mathrsfs,amssymb,amsthm,mathptm}
\usepackage{amsmath,amscd}
\usepackage{mathptm,pslatex}
\usepackage{color}
\usepackage[dvips]{graphicx}
\usepackage[all]{xy}
\usepackage{graphicx}
\usepackage{fancyhdr}

\begin{document}
\sloppy
\newcommand{\dickebox}{{\vrule height5pt width5pt depth0pt}}
\newtheorem{Def}{Definition}[section]
\newtheorem{Bsp}{Example}[section]
\newtheorem{Prop}[Def]{Proposition}
\newtheorem{Theo}[Def]{Theorem}
\newtheorem{Rem}[Def]{Remark}
\newtheorem{Lem}[Def]{Lemma}
\newtheorem{Koro}[Def]{Corollary}
\newtheorem{Claim}[Def]{Claim}

\newcommand{\cpx}[1]{#1^{\bullet}}
\newcommand{\lra}{\longrightarrow}
\newcommand{\lraf}[1]{\stackrel{#1}{\lra}}
 \newcommand{\ra}{\rightarrow}
\newcommand{\F}{\mathcal {F}}
\newcommand{\Hom}{{\rm Hom}}
\newcommand{\End}{{\rm End}}
\newcommand{\Ext}{{\rm Ext}}
\newcommand{\Tor}{{\rm Tor}}
\newcommand{\pd}{{\rm proj.dim}}
\newcommand{\inj}{{\rm inj}}
\newcommand{\lgd}{{l.{\rm gl.dim}}}
\newcommand{\gld}{{\rm gl.dim}}
\newcommand{\fd}{{\rm fin.dim}}
\newcommand{\Fd}{{\rm Fin.dim}}
\newcommand{\lfd}{l.{\rm Fin.dim}}
\newcommand{\rfd}{r.{{\rm Fin.dim}}}
\newcommand{\Mod}{{\rm Mod}}
\newcommand{\opp}{^{\rm op}}
\newcommand{\Proj}{{\rm Proj}}
\newcommand{\modcat}[1]{#1\mbox{{\rm -mod}}}
\newcommand{\pmodcat}[1]{#1\mbox{{\rm -proj}}}
\newcommand{\Pmodcat}[1]{#1\mbox{{\rm -Proj}}}
\newcommand{\injmodcat}[1]{#1\mbox{{\rm -inj}}}
\newcommand{\E}{{\rm E}_{\mathcal {F}}^{{\rm F},\Phi}}
\newcommand{\X}{ \mathscr{X}_{\mathcal {F}}^{{\rm F},\Phi}}
\newcommand{\Y}{\mathscr{Y}_{\mathcal {F}}^{{\rm F},\Phi}}
\newcommand{\A}{\mathcal {A}}
\newcommand{\C}{\mathscr{C}}
\newcommand{\K}{\mathscr{K}}
\newcommand{\D}{\mathscr{D}}
\newcommand{\Db}[1]{ \mathscr{D}^{\rm b}(#1)}
\newcommand{\Kb}[1]{{\mathscr K}^b(#1)}
\newcommand{\otimesL}{\otimes^{\rm\bf L}}
\newcommand{\otimesP}{\otimes^{\bullet}}
\newcommand{\rHom}{{\rm\bf R}{\rm Hom}}
\newcommand{\projdim}{\pd}
\newcommand{\stmodcat}[1]{#1\mbox{{\rm -{\underline{mod}}}}}
\newcommand{\Modcat}[1]{#1\mbox{{\rm -Mod}}}
\newcommand{\procat}[1]{#1\mbox{{\rm -proj}}}
\newcommand{\Tr}{{\rm Tr}}
\newcommand{\add}{{\rm add}}
\newcommand{\I}{{\rm Im}}
\newcommand{\Ker}{{\rm Ker}}
\newcommand{\EA}{\rm E^\Phi_\mathcal {A}}
\newcommand{\pro}{{\rm pro}}
\newcommand{\Coker}{{\rm Coker}}
\newcommand{\id}{{\rm id}}
\renewcommand{\labelenumi}{\Alph{enumi}}
\newcommand{\M}{\mathcal {M}}
\newcommand{\Mf}{\rm \mathcal {M}^f}
\newcommand{\rad}{{\rm rad}}
\newcommand{\injdim}{{\rm inj.dim}}

{\Large \bf
\begin{center}
Higher algebraic $K$-groups and $\mathcal D$-split sequences
\end{center}}
\medskip

\centerline{{\bf  Changchang Xi}}
\begin{center} School of Mathematical Sciences, Beijing Normal University, \\
Laboratory of Mathematics and Complex Systems, \\
100875 Beijing, People's Republic of  China \\ E-mail: xicc@bnu.edu.cn \\
\end{center}

\renewcommand{\thefootnote}{\alph{footnote}}
\setcounter{footnote}{-1} \footnote{2010 Mathematics Subject
Classification: Primary 19D50, 18E30; Secondary 13D15, 16S50.}
\renewcommand{\thefootnote}{\alph{footnote}}
\setcounter{footnote}{-1} \footnote{Keywords: Algebraic $K$-theory;
Derived equivalence; $\mathcal D$-split sequence; GV-ideal;
Mayer-Vietoris sequence.}

\begin{abstract}
In this paper, we use $\mathcal D$-split sequences and derived
equivalences to provide formulas for calculation of higher algebraic
$K$-groups (or mod-$p$ $K$-groups) of certain matrix subrings which
cover tiled orders, rings related to chains of Glaz-Vasconcelos
ideals, and some other classes of rings. In our results, we do not
assume any homological requirements on rings and ideals under
investigation, and therefore extend sharply many existing results of
this type in the algebraic $K$-theory literature to a more general
context.
\end{abstract}

\section{Introduction}
One of the fundamental questions in the algebraic $K$-theory of
rings is to understand and calculate higher algebraic $K$-groups
$K_n$ of rings, which were deeply developed in a very general
context by Quillen in \cite{Quillen} for exact categories and by
Waldhausen in \cite{wald} for Waldhausen categories. On the one
hand, the usual methods for computing $K_n$ may be the fundamental
theorem, splitting morphisms, or certain long exact sequences of
$K_n$-groups, namely, Mayer-Vietoris sequences, localization
sequences or excision. In this direction there is a lot of
literature (for example, see \cite{gkuk, kuk, swan1, weibel1,
weibel2, weibel}, and others). On the other hand, we know that
derived-equivalent rings share many common homological and numerical
features, in particular, they have the isomorphic higher algebraic
$K_n$-groups for all $n\ge 0$ (see \cite{DS}). This means that, in
order to understand the higher $K$-groups $K_n$ of a ring, one might
refer to another ring which is derived-equivalent to the given one,
and which may hopefully have a simple form so that its $K_n$-groups
can be determined easily. This idea, however, seems not much to be
benefited in the study of higher algebraic $K$-theory of rings,
especially in dealing with calculation of $K_n$-groups.

In the present note, we shall use ring extensions and derived
equivalences as reduction techniques to investigate the higher
algebraic $K_n$-groups of certain matrix subrings which include many
maximal orders, hereditary orders, tiled orders, endomorphism rings
of chains of Glaz-Vasconcelos ideals, and other classes of rings. To
produce such derived equivalences, we shall employ $\mathcal
D$-split sequences defined in \cite{hx2}. In this way, we reduce our
calculation inductively to that of certain triangular matrix rings.
The advantage of our method is: We not only drop all homological
conditions on rings and ideals under investigation, but also extend
many existing results (see \cite{bk, gkuk, keating}) of this type in
the literature to a more general context. Our main results in this
note can be stated as follows.

\begin{Theo} Suppose that $p\ge 2$, $m$ and $s$ are positive integers
such that $s$ is divisible by $p$. Let $R$ be a ${\mathbb
Z}/p^m\mathbb{Z}$-algebra with identity, and let $I, I_i$ and
$I_{ij}$ be (not necessarily projective) ideals of $R$. We denote by
$K_*(R)$ the $*$-th algebraic $K$-group of $R$ with $*\in
\mathbb{N}$.

$(1)$ If $I_{ij}\subseteq I$ for all $i,j$, $I_{kj}\subseteq I_{ij}$
for $k\le i$, $I_{ki}\subseteq I_{kj}$ for $j\le i$ and
$I_{ik}I_{kj}\subseteq I_{ij}$ for $i<k<j$, then
$$S:=\begin{pmatrix}
R & I_{12} &I_{13} & \cdots & I_{1\;n}\\
I & R & I_{23} & \cdots & I_{2\;n} \\
\vdots & \ddots &\ddots &\ddots  & \vdots\\
I & \cdots & I  & R& I_{n-1\;n}\\
I& \cdots &I & I & R
\end{pmatrix}$$ is a ring, and
$$  K_*(S)\otimes_{\mathbb Z}{\mathbb Z}[\frac{1}{s}]\simeq K_*(R)\otimes_{\mathbb Z}{\mathbb Z}[\frac{1}{s}]
\oplus \bigoplus_{j=2}^nK_*(R/I_{j-1\; j})\otimes_{\mathbb
Z}{\mathbb Z}[\frac{1}{s}].$$

$(2)$ For $2\le i\le n$, suppose that $R_i$ is a subring of $R$ with
the same identity. If $I_{i+1}\subseteq I_i\subseteq R_i$ for all
$i$, $I_j\subseteq I_{ij}\subseteq I$ for all $i,j$, and
$I_{ik}I_{kj}\subseteq I_{ij}$ for $j<k<i$, then
$$ T:=\begin{pmatrix}
R & I_2 &I_3 & \cdots  & I_n\\
I & R_2 & I_3 & \cdots & I_n\\
I & I_{32}& \ddots & \ddots & \vdots\\
\vdots & \vdots &\ddots & R_{n-1}& I_n\\
I & I_{n2} & \cdots & I_{n\, n-1} & R_n
 \end{pmatrix}$$ is a ring, and
$$  K_*(T)\otimes_{\mathbb Z}{\mathbb Z}[\frac{1}{s}]\simeq K_*(R)\otimes_{\mathbb Z}{\mathbb Z}[\frac{1}{s}]
\oplus \bigoplus_{j=2}^nK_*(R_j/I_j)\otimes_{\mathbb Z}{\mathbb
Z}[\frac{1}{s}].$$ \label{thm1}
\end{Theo}

As pointed out in Section \ref{modp} below, Theorem \ref{thm1} holds
true for the mod-$p$ $K$-groups $K_*(-,\mathbb{Z}/p\mathbb{Z})$ if
we assume in Theorem \ref{thm1} that $R$ is a
$\mathbb{Z}[\frac{1}{p}]$-algebra and $p\not\equiv 2 \; \mbox{(mod
4)}$, that is, under these two assumptions, one can replace
$K_*(-)\otimes_{\mathbb Z}\mathbb{Z}[\frac{1}{s}]$ by
$K_*(-,\mathbb{Z}/p\mathbb{Z})$ in Theorem \ref{thm1}.

The proof of the above result is based on the following observation.
Note that the assumptions in Theorem \ref{thm2}(2) below is weaker
than the ones in Theorem \ref{thm1}(2) above.

\begin{Theo} Let $R$ be a ring with identity, and let $I_{ij}$ be
(not necessarily projective) ideals of $R$. We denote by $K_*(R)$
the $*$-th algebraic $K$-group of $R$ with $*\in \mathbb{N}$.

$(1)$ If $I_{kj}\subseteq I_{ij}$ for $k\le i$, $I_{ki}\subseteq
I_{kj}$ for $j\le i$ and $I_{ik}I_{kj}\subseteq I_{ij}$ for $i<k<j$,
then
$$S:=\begin{pmatrix}
R & I_{12} &I_{13} & \cdots & I_{1\;n}\\
R & R & I_{23} & \cdots & I_{2\;n} \\
\vdots & \ddots &\ddots &\ddots  & \vdots\\
R & \cdots & R  & R& I_{n-1\;n}\\
R& \cdots &R & R & R
\end{pmatrix}$$ is a ring, and
$$  K_*(S)\simeq K_*(R)\oplus \bigoplus_{j=2}^nK_*(R/I_{j-1\; j}).$$

$(2)$  For $2\le i\le n$, suppose that $R_i$ is a subring of $R$
with the same identity, that $I_i \subseteq R_i$ is a right ideal of
$R_i$, and that $I_i$ is a left ideal of $R$. If $I_{i+1}\subseteq
I_i$ for all $i$, $I_j\subseteq I_{ij}$ for all $i,j$, and
$I_{ik}I_{kj}\subseteq I_{ij}$ for $j<k<i$, then
$$ T:=\begin{pmatrix}
R & I_2 &I_3  & \cdots  & I_n\\
R & R_2 & I_3 & \cdots  & I_n\\
R & I_{32}& \ddots & \ddots& \vdots\\
\vdots & \vdots &\ddots & R_{n-1}& I_n\\
R & I_{n2} & \cdots & I_{n\, n-1} & R_n
\end{pmatrix}$$ is a ring, and
$$  K_*(T)\simeq K_*(R)\oplus \bigoplus_{j=2}^nK_*(R_j/I_j).$$
\label{thm2}
\end{Theo}

The strategy  of our proofs of the theorems is first to use ring
extensions, which are motivated from \cite{xi2}, and then to combine
$K$-groups in Mayer-Vietoris sequences with $K$-groups of rings
which are linked by derived equivalences produced from certain
$\mathcal D$-split sequences.

This note is organized as follows. In Section \ref{sect2}, we recall
some definitions and elementary facts on derived equivalences needed
in the later proofs. In Section \ref{extension}, we construct
$\mathcal D$-split sequences by ring extensions and calculate the
endomorphism rings of tilting modules related to these sequences. In
Section \ref{ktheory}, we prove the main results and state some of
its consequences. Our proofs of the above results also give an
explanation of the multiplicity factor $n-1$ in the isomorphisms of
$K_n$-groups of the rings in \cite{gkuk} and \cite{keating}. In
Section \ref{lowerK}, we calculate $K_0$ and $K_1$ for some matrix
subrings which are not covered by the main results. In Section
\ref{modp}, we show that the main result Theorem \ref{thm1} holds
for mod-$p$ $K$-theory by outlining the key ingredients of its
proof. In Section \ref{example}, we give some examples to show how
our method can work, here GV-ideals in commutative rings enter into
our play. These examples demonstrate also that the matrix rings
studied in Section \ref{extension} really occur, as the endomorphism
rings of  chains of GV-ideals, in the field of commutative algebra.

\section{Preliminaries\label{sect2}}

Let $A$ be a ring with identity. By an $A$-module we mean a left
$A$-module. Let $A$-Mod (respectively, $A$-mod) denote the category
of all (respectively, finitely generated) left $A$-modules.
Similarly, by $A$-Proj (respectively, $A$-proj) we denote the full
subcategory of all (finitely generated) projective $A$-modules in
$R$-Mod. For an $A$-module $M$, we denote by $\pd(_AM)$ the
projective dimension of $M$. Let $\Kb{\pmodcat{A}}$ be the homotopy
category of the additive category $A$-proj. The unbounded derived
category of $A$-Mod is denoted by $\D(A)$, whereas the bounded
derived category of $A$-Mod is denoted by $\Db{A}$. We say that two
rings $A$ and $B$ are derived-equivalent if $\D(A)$ and $\D(B)$ are
derived-equivalent as triangulated categories. It is well-known that
if $\Db{A}$ and $\Db{B}$ are derived-equivalent as triangulated
categories then $\D(A)$ and $\D(B)$ are derived-equivalent as
triangulated categories.

Given a additive category $\mathcal C$ and an object $X$ in
$\mathcal C$, we denote by $\add(X)$ the full subcategory of
$\mathcal C$ consisting of all objects which are direct summands of
direct sums of finitely many copies of $X$.

For derived equivalences, Rickard's Morita theory \cite{rickard} is
very useful.

\begin{Theo}{\rm\cite{rickard}}
For two rings $A$ and $B$ with identity, the following are
equivalent:

$(a)$ $\Db{A}$ and $\Db{B}$ are equivalent as triangulated
categories.

$(b)$ $\Kb{\pmodcat{A}}$ and $\Kb{\pmodcat{B}}$ are equivalent as
triangulated categories.

$(c)$ $B\simeq \End_{\Kb{\pmodcat{A}}}(\cpx{T})$, where $\cpx{T}$ is
a complex in $\Kb{\pmodcat{A}}$ satisfying

\parindent=1cm $(1)$  $\Hom(\cpx{T},\cpx{T}[i])=0$ for $i\neq 0$, and

$(2)$ $\add(\cpx{T})$ generates $\Kb{\pmodcat{A}}$ as a triangulated
category.\label{rickard}
\end{Theo}

For derived equivalences, it is shown in \cite{DS} that the
algebraic $K$-theory is an invariant. Recall that, for a ring $A$
with identity, $K_n(A)$ denotes the $n$-th homotopy group  of a
certain space $K(A)$ produced by one¡¯s favorite $K$-theory defined
for each $n\in \mathbb{N}$ (see \cite{Quillen}, \cite{wald},
\cite{weibel}).

\begin{Theo}{\rm \cite{DS}} If two rings $A$ and $B$ with identity are derived-equivalent, then their algebraic
$K$-groups are isomorphic: $K_*(A)\simeq K_*(B)$ for all $*\in
\mathbb{N}$. \label{dg}
\end{Theo}

As is known, Morita equivalences are derived equivalences. Thus, if
$A$ and $B$ are Morita equivalent, then their algebraic $K$-groups
are isomorphic.

Another special class of derived equivalences can be constructed by
tilting modules initialled from the representation theory of
finite-dimensional algebras (for example, see \cite{BB}). Recall
that a module $T$ over a ring $A$ is called a tilting module if the
following three conditions are satisfied:

$(1)$ $T$ has a finite projective resolution $0\lra
P_n\lra\cdots\lra P_0\lra T\lra 0$, where each $P_i$ is a finitely
generated projective $A$-module;

$(2)$ $\Ext^i_A(T,T)=0$ for all $i>0$, and

$(3)$ there is an exact sequence $0\lra A\lra T^0\lra\cdots\lra
T^m\lra 0$ of $A$-modules with each $T^i$ in $\add(T)$ which is a
full subcategory of $A$-mod.

\medskip
Note that, for a tilting module $T$, the projective resolution
$\cpx{P}$ of $T$ satisfies  (1) and (2) of Theorem \ref{rickard}(c).
Thus, if $_AT$ is a tilting $A$-module then $A$ and $\End_A(T)$ are
derived-equivalent. To produce tilting modules, one may use the
notion of $\mathcal D$-split sequences. Now let us recall the
definition of $\mathcal D$-split sequences from \cite{hx2}.

Let $\cal C$ be  an additive category and $\cal D$ a full
subcategory of $\cal C$. A sequence
$$X\stackrel{f}{\longrightarrow}M\stackrel{g}{\longrightarrow}Y$$
of morphisms between objects in $\mathcal C$ is called a $\cal
D$-split sequence if

$(1)$ $M\in {\cal D}$,

$(2)$ $f$ is a left $\cal D$-approximation of $X$, that is,
$\Hom_{\mathcal C}(f,D'): \Hom_{\mathcal C}(M,D')\lra\Hom_{\mathcal
C}(X,D')$ is surjective for all $D'\in{\mathcal D}$, and $g$ is a
right $\cal D$-approximation of $Y$, that is, $\Hom_{\mathcal
C}(D',g): \Hom_{\mathcal C}(D',M)\lra\Hom_{\mathcal C}(D',Y)$ is
surjective for all $D'\in \mathcal{D}$, and

$(3)$ $f$ is a kernel of $g$, and $g$ is a cokernel of $f$.

\medskip
Examples of $\mathcal D$-split sequences include Auslander-Reiten
sequences and short exact sequences of the form $0\ra X\ra P\ra Y\ra
0$ in $A$-Mod with $P$ projective-injective. A non-example is the
sequence $0\ra \mathbb{Z}\ra \mathbb{Q}\ra \mathbb{Q}/\mathbb{Z}\ra
0$ of ableilan groups, that is, this sequence is not an
$\add(\mathbb{_{\mathbb{Z}}Q})$-split sequence. For more examples,
one may find in \cite{hx2}, and also in the next section as well as
in the last section of the present paper.

Given a $\cal D$-split sequence $X\ra M'\ra Y$, with $\cal D$ =
$\add(M)$ for $M$ an object in $\mathcal C$, it is shown in
\cite{hx2} that there is a tilting module $T$ over $\End_{\mathcal
C}(X\oplus M)$ of projective dimension at most $1$ such that
$\End(T)$ is isomorphic to $\End_{\mathcal C}(M\oplus Y)$. Thus
$\End_{\mathcal C}(X\oplus M)$ and $\End_{\mathcal C}(M\oplus Y)$
are derived-equivalent, and have the isomorphic algebraic $K$-theory
by Theorem \ref{dg}.

\section{Ring extensions and derived equivalences\label{extension}}

Ring extensions were used in \cite{xi2} to study the finitistic
dimensions of algebras. In this section, we shall use ring
extensions to construct $\mathcal D$-split sequences which will be
applied to calculation of the algebraic $K$-groups of rings in the
next section.

We first establish the following  general fact.

\begin{Lem}  Let $B\subseteq A$ be an extension of rings with the same identity.

$(1)$ If $\Ext^1_B(_BA, {}_BB)=0$, then the sequence
$$(*)\qquad 0\lra B\lra A\lra A/B\lra 0$$ is an $\add(_BA)$-split sequence in
$B\Modcat$. Thus $\End_B(_BB\oplus {}_BA)$ and $\End_B(_BA\oplus
A/B)$ are derived-equivalent.

$(2)$ If $_BA$ is projective, then the above sequence is an
$\add(_BA)$-split sequence.

$(3)$ Suppose that $\Ext^1_B(_BA, {}_BA)=0$. If $_BA$ is finitely
presented with $\pd(_BA)\le 1$ (for instance, $_BA$ is projective
and finitely generated), then $A\oplus A/B$ is a tilting $B$-module
of projective dimension at most $1$. In particular, $\End_B(A\oplus
A/B)$ is derived-equivalent to $B$.
 \label{lem1}
\end{Lem}

{\it Proof.} (1) We have the following exact sequence
$$ 0\ra \Hom_B(_BA,B)\lra \Hom_B(A,A)\lra\Hom_B(A, A/B)\lra \Ext^1_B(A,B)\lra\Ext^1_B(A,A).$$
The condition $\Ext^1_B(A,B)=0$ implies that the canonical
surjection $A\ra A/B$ is a right $\add(_BA)$-approximation of $A/B$.
To see that the inclusion $B\lra A$ is a left
$\add(_BA)$-approximation of $B$, we note that each homomorphism
from $_BB$ to $_BA$ is given by an element $a$ in $A$. Thus it can
be extended to a homomorphism from $_BA$ to $_BA$ by the right
multiplication of $a$. Clearly, one can check that this is also true
for any homomorphism from $_BB$ to a direct summands of $_BA$. Thus
we see that the inclusion map from $B$ to $A$ is a left
$\add(_BA)$-approximation of $B$. Thus ($*$) is an $\add(_BA)$-split
sequence in $B\Modcat$, and therefore $\End_B(_BB\oplus {}_BA)$ and
$\End_B(_BA\oplus A/B)$ are derived-equivalent by \cite[Theorem
1.1]{hx2}. This finishes the proof of Lemma \ref{lem1}(1).

(2) is a special case of (1).

(3) Let $T:= {}_BA\oplus \; A/B$. Since $_BA$ is finitely presented
of projective dimension at most one, there is an exact sequence
$0\ra P_1\ra P_0\ra {}_BA\ra 0$ such that $P_i$ are finitely
generated projective $B$-modules and the following diagram is
commutative: {\footnotesize $$ \begin{CD}@.0 @. 0\\@.@VVV @VVV\\
@. P_1@= P_1 \\ @. @VVV @VVV\\
0@>>>P@>>> P_0@>>>A/B@>>>0\\
@.@VVV @VVV @| \\
0@>>> B@>>>A@>>>A/B@>>> 0.\\
@. @VVV @VVV \\ @.0 @.0
\end{CD} $$} From this diagram we see that $T$ is
finitely presented of projective dimension at most one. Thus the
conditions (1) and (3) in the definition of tilting modules are
satisfied. It remains to show $\Ext^1_B(A\oplus A/B, A\oplus
A/B)=0$. This is equivalent to that $\Ext^1_B(A, A/B)=0$,
$\Ext^1_B(A/B, A/B)=0$  and $\Ext^1_B(A/B,A)=0$ since
$\Ext^1_B(A,A)=0$ by assumption.

Indeed, we have seen that the inclusion map $\lambda$ from $B$ into
$A$ is always a left $\add(_BA)$-approximation of $_BB$. Thus the
induced map $\lambda^*:=\Hom_{B}(\lambda, A)$ is surjective. Hence,
by applying $\Hom_B(-,A)$ to the canonical exact sequence ($*$), we
get an exact sequence
$$ 0\ra \Hom_B(A/B,A)\lra \Hom_B(A,A)\lraf{\lambda^*}\Hom_B(B,A)\lra \Ext^1_B(A/B,A)\lra \Ext^1_B(A,A),$$
which shows $\Ext^1_B(A/B,A)=0$. If we apply $\Hom_B(A/B,-)$ to the
canonical exact sequence, then we get an exact sequence:
$$ \Ext^1_B(A/B, B)\lra \Ext^1_B(A/B,A)\lra \Ext^1_B(A/B, A/B)\lra 0$$ since the projective dimension of $A/B$ is at most $1$.
This implies $\Ext^1_B(A/B,A/B)=0$. Similarly, applying
$\Hom_B(A,-)$ to the canonical exact sequence ($*$), we can deduce
$\Ext^1_B(A,A/B)=0$. Thus we complete the proof of (3). $\square$

{\it Remark}. Sometimes the following observation is useful for
getting $\mathcal D$-split sequences: Suppose that $e$ and $f$ are
idempotent elements in a ring $R$ and $a\in eRf$. Then the right
multiplication map $Re\ra Rf$, defined by $x\mapsto xa$ for $x\in
Re$, is a left $\add(Rf)$-approximation of $Re$ if and only if
$eRf=afRf$. Thus, if the right multiplication map is injective, then
$0\ra Re\ra Rf\ra Rf/Rea\ra 0$ is an $\add(Rf)$-split sequence if
and only if $eRf=afRf$. For instance, the sequence  $0\ra
\mathbb{Z}\lraf{\cdot 2}\mathbb{Z}\ra {\mathbb Z}/2{\mathbb Z}\ra 0$
is not an $\add(\mathbb Z)$-split sequence.

\medskip Let us mention an example of ring extensions which satisfy the
conditions in Lemma \ref{lem1}. Recall that an extension $B\subseteq
A$ of rings is called a quasi-Frobenius extension if $_BA$ is
finitely generated and projective, and the bimodule $_AA_B$ is
isomorphic to a direct summand of the direct sum of finitely many
copies of $_A\Hom_B(_BA_A, {}_BB)_B$. Thus each quasi-Frobenius
extension $B\subseteq A$ provides an $\add(_BA)$-split sequence
$$ 0\lra B\lra A\lra A/B\lra 0,$$ and a tilting $B$-module $A\oplus A/B$ by Lemma \ref{lem1}.

\medskip
Now we consider some consequences of Lemma \ref{lem1}, which are
needed in the next section.

Let $R$ be a ring with identity and $I_{ij}$ ideals in $R$ with
$1\le i<j\le n$, such that

(1) $I_{kj}\subseteq I_{ij}$ for $k\le i$,

(2) $I_{ki}\subseteq I_{kj}$ for $j\le i$, and

(3) $I_{ik}I_{kj}\subseteq I_{ij}$ for $i < k <j$. Then
$$B:=\begin{pmatrix}
R & I_{12} &I_{13} & \cdots & I_{1\;n}\\
R & R & I_{23} & \cdots & I_{2\;n} \\
\vdots & \vdots &\ddots &\ddots  & \vdots\\
R & R &  \cdots & R& I_{n-1\;n}\\
R& R& \cdots& R & R
\end{pmatrix}$$
is a ring. The rings of this form include tiled triangular orders
and maximal orders \cite{cr}. They occur also as the endomorphism
rings of chains of Glaz-Vasconcelos ideals of rings, see Section
\ref{example}.

The following lemma shows that we may use derived equivalences to
simplify the ring $B$.

\begin{Lem} Let $B$ be the ring defined above. Then $B$ is
derived-equivalent to
$$C:=\begin{pmatrix}
R & I_{12} &I_{13} & \cdots & I_{1\, n-1} & I_{1\, n-1}/I_{1n}\\
R & R & I_{23} & \cdots & I_{2\, n-1} & I_{2\,n-1}/I_{2n}\\
R & R & R &\ddots & \vdots & \vdots\\
\vdots & \vdots &\vdots &\ddots  & I_{n-2\, n-1}& I_{n-2\, n-1}/I_{n-2\, n}\\
R & R & R&\cdots & R & R/I_{n-1\;n}\\
0 & 0 & 0& \cdots & 0& R/I_{n-1\;n}
\end{pmatrix}.$$
\label{tde}
\end{Lem}

{\it Proof.} We make the following conventions on notations. Let
$S=M_n(R)$, the $n \times n$ matrix ring over $R$. Let $e_i$ be the
$n\times n$ matrix with $1_R$ in $(i,i)$-entry and zero in other
entries. For convenience, we denote by $e_{i,j}(x)$ the matrix with
$x$ in $(i,j)$-position, and zero in other positions, and by
$B_{ij}$ the $(i,j)$-component of the matrix subring $B$ of $S$,
that is, the set of  $(i,j)$-entries of all matrices in $B$. We
define
$$A:=\begin{pmatrix}
  R        &  I_{12} &  \cdots     & I_{1\;n-1}      &  I_{1\;n-1}\\
  R        &  R          &  \ddots     &  \vdots         &  \vdots\\
  \vdots   & \vdots      & \ddots      &  I_{n-2\; n-1}&  I_{n-2\; n-1}\\
  R        & R           & \cdots      & R               & R \\
  R        & R           &  \cdots           & R               & R
\end{pmatrix}.$$
Note that the only difference between $A$ and $B$ is the last
column. We can verify that $A$ is a ring containing $B$ as a
subring.

Clearly, as a left $B$-module, $_BA\simeq Be_1\oplus \cdots \oplus
Be_{n-1}\oplus Be_{n-1}$. Thus $_BA$ is finitely generated and
projective. Furthermore, it follows that $B$ is Morita equivalent to
$\End_B(B\oplus {}_BA)$ and that the latter is derived-equivalent to
$\End_B(_BA\oplus A/B)$ by Lemma \ref{lem1}. Thus $B$ is
derived-equivalent to $\End_B(Be_1\oplus \cdots\oplus Be_{n-1}\oplus
\; Ae_n/Be_n)$. For simplicity, we denote by $Q$ the $B$-module
$Ae_n/Be_n$. Note that $Ae_n\simeq Be_{n-1}$ as $B$-modules, and
that we have a canonical exact sequence: $$0\lra
Be_n\lraf{\lambda}Be_{n-1}\lraf{\pi} Q\lra 0,$$ where $\lambda$ is
the composition of the inclusion of $Be_{n}$ into $Ae_n$ with the
right multiplication $\cdot e_{n,n-1}$, and $\pi$ is the composition
of the right multiplication of $\cdot e_{n-1,n}$ with the canonical
surjective map $Ae_n\ra Ae_n/Be_n$.

In the following, we shall prove that $\End_B(Be_1\oplus
\cdots\oplus Be_{n-1}\oplus Q)$ is isomorphic to $C$.

First, we define a map $\varphi: R\lra \Hom_B(Q,Q)$ as follows: For
$b\in R$, let $\cdot e_nbe_n$ be the right multiplication map from
$Be_n$ to $Be_n$. This is well-defined by our assumptions. Also, let
$\cdot e_{n-1}be_{n-1}$ be the right multiplication map from
$Be_{n-1}$ to itself. Then we see that $\lambda(\cdot
e_{n-1}be_{n-1})=(\cdot e_nbe_n)\lambda$. So, there is a unique
$\alpha\in\Hom_B(Q,Q)$ making the following diagram commutative:

$$\begin{CD}
@.0@>>>Be_n@>{\lambda}>>Be_{n-1}@>{\pi}>> Q@>>>0\\
(**) @. @. @V{\cdot e_nbe_{n}}VV @V{\cdot e_{n-1}be_{n-1}}VV @VV{\alpha}V \\
@. 0@>>>Be_n@>{\lambda}>>Be_{n-1}@>{\pi_j}>> Q@>>>0.\\
\end{CD}$$
Hence, we can define the image of $b$ under $\varphi$ is $\alpha$.
Clearly, if $b,b'\in R$, then $(b+b')\varphi=(b)\varphi
+(b')\varphi$. Since $e_n(bb')e_n=e_nbe_nb'e_n$, we also have
$(bb')\varphi= (b\varphi)(b'\varphi)$. Thus $\varphi$ is a
homomorphism of rings.

Now, we calculate the kernel of $\varphi$. Suppose $b\in R$ such
that $\alpha = b\varphi=0$. Then the map $\cdot e_{n-1}be_{n-1}$
factorizes through $\lambda$. This means that there is an element
$r\in B_{n-1\,n}$ such that $\cdot e_{n-1}be_{n-1}=(\cdot
e_{n-1}re_n)\lambda$ and $\cdot e_nbe_n = \lambda(\cdot
e_{n-1}re_n)$. Hence $b=r\in B_{n-1\, n}.$ Thus
Ker$(\varphi)\subseteq B_{n n}\cap B_{n-1\,n}=I_{n-1\,n}$. Since any
map $\cdot e_nbe_n$ from $Be_n$ to $Be_n$ with $b\in B_{n n}\cap
B_{n-1 n}$ factorizes through $\lambda$, the corresponding $\alpha$
is zero. Hence Ker$(\varphi)$ is $B_{n-1\; n}$.

Given an element $\alpha\in\Hom_B(Q,Q)$, we may form the following
commutative diagram in $B$-Mod:
$$\begin{CD}
0@>>>Be_n@>{\lambda}>>Be_{n-1}@>{\pi}>> Q@>>>0\\
@. @V{b}VV @V{a}VV @VV{\alpha}V \\
0@>>>Be_n@>{\lambda}>>Be_{n-1}@>{\pi}>> Q@>>>0.\\
\end{CD}$$
Note that the homomorphism $a$ exists and makes the right square of
the above diagram commutative. Thus we have a homomorphism $b$
making the left square commutative. We may identify $a$ with an
element in $B_{n-1\,n-1}$, say $a=\cdot e_{n-1,n-1}(r)$ with $r\in
B_{n-1\,n-1}$, and identify $b$ with an element in $B_{nn}$, say
$b=\cdot e_{n,n}(s)$ with $s\in B_{nn}$. The commutativity of the
left square means that $r=s\in B_{nn}$. This means that $\varphi$ is
surjective. Thus $\End_B(Q)\simeq R/I_{n-1\,n}.$

If we apply $\Hom_B(-,Be_{j})$ to $(*)$ for $1\le j\le n-1$ and use
Lemma \ref{lem1}(3), we have the following exact commutative diagram
with $e_{n, n-1}\cdot $ an isomorphism: {\footnotesize $$\begin{CD}
0@>>>\Hom_B(Q,Be_{j})@>{(\pi)_*}>>\Hom_B(Be_{n-1},Be_{j})@>{(\lambda)_*}>> \Hom_B(Be_n,Be_j)@>>>0\\
@. @. @V{\simeq}VV  @VV{\simeq}V \\
@.@. e_{n-1}Be_{j}@>{e_{n,n-1}\cdot}>> e_nBe_j.\\
\end{CD}$$}
Thus $\Hom_B(Q,Be_j)=0$ for all $1\le j\le n-1$.

If we apply $\Hom_B(Be_j,-)$ to the exact sequence ($*$) for $1\le
j\le n-1$, we get an exact sequence

$$ 0\lra \Hom_B(Be_j, Be_n)\lra \Hom_B(Be_j,B_{n-1})\lra \Hom_B(Be_j,
Q)\lra 0,$$ which shows that $\Hom_B(Be_j,Q)\simeq
B_{j\,n-1}/B_{j\,n}=B_{j\,n-1}/I_{j\; n-1}.$

Now we identify $\Hom_B(Be_j,Be_i)$ with $e_jBe_i$ for all $1\le
i,j\le n-1$, and $\Hom_B(Be_j,Q)$ with
$B_{j\,n-1}/B_{j\,n}=B_{j\,n-1}/I_{j\; n-1}.$ Then we can see that
$\End_B(Be_1\oplus \cdots\oplus Be_{n-1}\oplus Q)$ is isomorphic to
$C$. This finishes the proof of Lemma \ref{tde}. $\square$

A special case of Lemma \ref{tde} is the ring considered in
\cite{gkuk} under certain homological assumptions and finiteness
condition. Here we start with a more general setting and drop all
homological conditions on ideals as well as finiteness condition on
quotients.

Let $R$ be a ring with identity, and $I$ an arbitrary ideal in $R$.
We consider the ring of the following form

$$B:=\begin{pmatrix}
  R        &  I^{t_{12}} &  \cdots     & I^{t_{1n}}    \\
  R        &  R          &  \ddots     &  \vdots  \\
  \vdots   & \vdots      & \ddots      &  I^{t_{n-1\; n}}   \\
  R        & R           & \cdots      & R
\end{pmatrix},$$ where $t_{ij}$ are positive integers. Note that the conditions for $B$ to be a ring are

(1) $t_{ij}\le t_{i\; j+1}, \; t_{i+1\; j}\le t_{ij}$ for $i<j$, and

(2) $t_{ij}\le t_{ik}+t_{kj}$ for $i<k<j$.

\medskip
The next result follows immediately from Lemma \ref{tde}.

\begin{Lem} Assume that the above defined $B$ is a ring. Then $B$ is
derived-equivalent to
$$C:=\begin{pmatrix}
  R        &  I^{t_{12}} &  \cdots     & I^{t_{1\;n-1}}     & I^{t_{1\;n-1}}/I^{t_{1\; n}}\\
  R        &  R          &  \ddots     &  \vdots            & \vdots \\
  \vdots   & \vdots      & \ddots      &  I^{t_{n-2\; n-1}} &  I^{t_{n-2\; n-1}}/I^{t_{n-2\;n}}   \\
  R        & R           & \cdots      &  R                 & R/I^{t_{n-1\; n}}\\
  0        & 0           &  0          &  0                 &
  R/I^{t_{n- 1\; n}}
\end{pmatrix}.$$
\label{lem2}
\end{Lem}

\medskip
Next, we consider a variation of the ring $B$ in Lemma \ref{tde},
which was considered in \cite{kk, chyp} and cover some tiled orders
in \cite{cr}, and many other cases, for example, rings in \cite{kk},
and some Auslander-regular, Cohen-Macaulay rings (not necessarily
maximal orders, see \cite{stafford}).

Let $R$ be a ring with identity. Suppose that $R_i$ is a subring of
$R$ with the same identity for $2 \leq i \leq n$, that $I_i$ is a
left ideal of $R$ for all $2\leq i\leq n$, and that $I_{ij}$ is
ideal of $R$, with $2\leq j<i\le n$, which satisfies the following
conditions:

(1) $I_i\subseteq R_i$ is a right ideal of $R_i$ for all $i$,

(2) $I_n\subseteq I_{n-1}\subseteq\cdots\subseteq I_2,$

(3) $ I_j\subseteq I_{ij}$ for all $i,j$, and

(4)  $I_{ik}I_{kj}\subseteq I_{ij}$ for $j<k<i$.

\noindent Here we do not assume that $I_i$ is projective as a left
$R$-module, nor that $I_i$ is an ideal of $R$. Nevertheless one can
check that

$$B :=
 \begin{pmatrix}
R & I_2 &I_3 & \cdots & I_{n-1} & I_n\\
R & R_2 & I_3 & \cdots & I_{n-1} & I_n\\
R & I_{32}& R_3& \ddots & \vdots & I_n\\
R & I_{42} & I_{43}& \ddots & I_{n-1} & \vdots\\
\vdots & \vdots &\vdots & \ddots & R_{n-1}& I_n\\
R & I_{n2} & I_{n3} & \cdots & I_{n\, n-1} & R_n
 \end{pmatrix}, \qquad
A :=
\begin{pmatrix}
R & R &I_3 & \cdots & I_{n-1} & I_n\\
R & R & I_3 & \cdots & I_{n-1} & I_n\\
R & R& R_3& \ddots & \vdots & I_n\\
\vdots & \vdots &\vdots &\ddots  & I_{n-1}& \vdots\\
R & R & I_{n-1\;3}&\cdots & R_{n-1} & I_n\\
R & R & I_{n 3}& \cdots & I_{n\,n-1} & R_n
\end{pmatrix},$$

$$C :=
\begin{pmatrix}
R_2/I_2      & 0      &   0         & \cdots  & 0          & 0\\
R/I_2        & R      &   I_3       & I_4     & \cdots     & I_n\\
R/I_{32}     & R      &   I_3       & I_4     & \cdots     & I_n\\
R/I_{42}     & R      & R_3         & I_4     & \cdots     & I_n\\
\vdots       & \vdots &\vdots       &\ddots   & \ddots     & \vdots\\
R/I_{n-1\,2} & R      & I_{n-1,3\;3}&\cdots   & R_{n-1}    & I_n\\
R/I_{n 2}    &  R     & I_{n3}      & \cdots  & I_{n\,n-1} & R_n
\end{pmatrix}$$with the usual matrix addition and multiplication form three rings
with identity. Note that only the second column of $A$ is different
from the one of $B$.

We define a $B$-module $Q$ as follows:
$$ 0\lra Be_2\lraf{\lambda} Be_1\lraf{\pi} Q\lra 0,$$
where $\lambda$ is a composition of the inclusion $Be_2\ra Ae_2$
with the isomorphism $Ae_2\simeq Be_1$ as $B$-modules and where
$\pi$ is the cokernel of $\lambda$.

Now, we consider the endomorphism ring $\End_B(Q\oplus Be_1\oplus
Be_3\oplus \cdots\oplus Be_n)$. By a proof similar to that of Lemma
\ref{tde}, one can show that the following lemma is true. We leave
the details of its proof to the reader.

\begin{Lem}  The above defined rings $B$
and $C$ are derived-equivalent. \label{de}
\end{Lem}

An alternative proof of Lemma \ref{de} can be found in \cite[Theorem
5.2]{chyp}, where $A$ is replaced by the $n\times n$ matrix ring
over $R$.

\section{Higher algebraic $K$-theory of matrix subrings\label{ktheory}}

In the algebraic $K$-theory of rings, the calculation of higher
algebraic $K$-groups $K_n$ seems to be one of the interesting and
hard problems. In this section, we shall provide formulas for
computation of the $K_n$-groups of certain rings by applying the
results in the previous section. Our computation is based the
philosophy that derived equivalences of rings preserve the
$K$-theory and $G$-theory (see \cite{DS}), thus one can transfer the
calculation of $K_n$ of a ring to that of another ring which is
derived-equivalent to and may be much more simpler than the original
one. In the literature, there are many papers dealing with
$K_n$-groups by exploiting excision, Mayer-Vietoris exact sequences
or other related sequences (for example, see \cite{gkuk},
\cite{kuk}, \cite{swan1}, \cite{weibel2}, \cite{weibel}). However,
it seems that there are few papers using derived equivalences to
calculate the higher algebraic $K$-groups. In the present section we
shall show that sometimes our philosophy works powerfully though it
may be difficult to find derived equivalences in general. For some
new advances in constructing derived equivalences, we refer the
reader to recent papers \cite{hkx, hx4}.

Let $R$ be a ring with identity. We denote by $K_*(R)$ the series of
algebraic $K$-groups of $R$ with $*\in\{0,1,2, \cdots, \}$. The
algebraic $K$-theory of matrix-like rings has been of interest since
a long time. In \cite{bk}, Berrick and Keating showed the following
result.

\begin{Lem}{\rm \cite{bk}} If $R_i$ is a ring with identity for $i=1,2$,
and if $M$ is an $R_1$-$R_2$-bimodule, then, for the triangular
matrix ring
$$ S=\begin{pmatrix}
  R_1             &  M    \\
  0    & R_2      &     \\
  \end{pmatrix},$$
there is an isomorphism of $K$-groups: $ K_n(S)\simeq K_n(R_1)\oplus
K_n(R_2)$ for all integers $n\in\mathbb Z$.  Moreover, this
isomorphism is induced from the canonical inclusion of $R_1\oplus
R_2$ into $S$.\label{k1}
\end{Lem}
For $n=0,$ this is classical. For $n = 1,2$, this was already shown
by Dennis and Geller in 1976. We remark that Lemma \ref{k1} can be
used to calculate the higher algebraic $K$-groups of algebras
associated to finite $EI$-categories, or more generally, of
``triangular" Artin algebras. Recall that an Artin algebra $A$ over
a commutative Artin ring is said to be triangular if the set of
non-isomorphic indecomposable projective $A$-modules can be ordered
as $P_1, P_2, \cdots, P_n$ such that $\Hom_A(P_j,P_i)=0$ for all
$j>i$. In this case, we have $K_*(A)\simeq
\bigoplus_{j+1}^nK_*(\End_A(P_j))$ by Lemma \ref{k1}.

For a matrix ring of the form

$$T={\begin{pmatrix}
  R        &  I &  \cdots     & I    \\
  R        &  R          &  \ddots     &  \vdots  \\
  \vdots   & \vdots      & \ddots      &  I   \\
  R        & R           & \cdots      & R
\end{pmatrix},}_{n\times n}$$ where $R$ is a ring and $I$ is an ideal in $R$ such that the $R$-modules
$_RI$ and $I_R$ are projective, it was shown by Keating in
\cite{keating} that there is an isomorphism of $K$-theory:
$$ K_*(T)\simeq K_*(R)\oplus (n-1)K_*(R/I).$$
In \cite{keating}, the author also considered the so-called trivial
extension of a ring by a bimodule. It was shown that if $T$ is the
trivial extension of a ring $R$ by an $R$-bimodule $M$, then
$K_*(T)\simeq K_*(R)$ provided that $M$ has finite projective
dimension as a left $T$-module. Here the condition on $M$ in this
statement is necessary. See the counterexample $T:=k[x]/(x^2)$ which
is the trivial extension of $k$ by $k$, where $k$ is any field.

Recently, as a kind of generalization of the above result of
Keating, the authors of \cite{gkuk} consider the following matrix
ring: Let $I$ be an ideal of a ${\mathbb Z}_p$-algebra $R$ with
identity, where $\mathbb{Z}_p$ is the $p$-adic integers (or,
equivalently,
$\mathbb{Z}_p=\displaystyle\varprojlim_{n}\mathbb{Z}/p^n\mathbb{Z}$),
and define

$$S=\begin{pmatrix}
  R        &  I^{t_{12}} &  \cdots     & I^{t_{1n}}    \\
  R        &  R          &  \ddots     &  \vdots  \\
  \vdots   & \vdots      & \ddots      &  I^{t_{n-1\; n}}   \\
  R        & R           & \cdots      & R
\end{pmatrix},$$ where $t_{ij}$ are positive integers. Assume that
$S$ is a ring and that $R/I^n$ is a finite ring for all $n\ge 1$. If
both $_RI$ and $I_R$ are projective, it is proved in \cite{gkuk}
that the following isomorphism of the algebraic $K$-theory holds:
$$  K_*(S)(1/s)\simeq K_*(R)(1/s)\oplus (n-1)K_*(R/I)(1/s),$$
where $s$ is any rational integer such that $p$ divides $s$, and
where $G(1/s)$ denotes the group $G\otimes_{\mathbb Z}{\mathbb
Z}[\frac{1}{s}]$ for an abelian group $G$.

We shall use our results in the previous section to extend all
results on matrix rings mentioned above without any homological
conditions on rings and ideals under investigation. Our proofs also
explain the reason why the multiplicity $n-1$ appears in the above
mentioned isomorphisms of the higher algebraic $K$-theory.

 \begin{Lem} Let $B$ be the matrix ring
 $$B:=\begin{pmatrix}
R & I_{12} &I_{13} & \cdots & I_{1\;n}\\
R & R & I_{23} & \cdots & I_{2\;n} \\
\vdots & \vdots &\ddots &\ddots  & \vdots\\
R & R &  \cdots & R& I_{n-1\;n}\\
R& R& \cdots& R & R
\end{pmatrix}$$
 defined in Lemma \ref{tde}. Then
 $$  K_*(B)\simeq K_*(R)\oplus \bigoplus_{j=1}^{n-1}K_*(R/I_{j\; j+1}).$$
\label{ktheory1}
\end{Lem}

{\it Proof.} We show this lemma by induction on $n$. By Theorem
\ref{dg} (see \cite{DS}), the algebraic $K$-theory and $G$-theory
are invariant under derived equivalences. So, by Lemma \ref{tde}, we
have $K_*(B)\simeq K_*(C)$ (for notation see Section
\ref{extension}). Now it follows from Lemma \ref{k1} that
$K_*(C)\simeq K_*(R/I_{n-1\, n})\oplus K_*(B_{n-1})$, where
$B_{n-1}$ is the $(n-1)\times (n-1)$ left upper corner matrix
subring of $B$. By induction, we have $K_*(B_{n-1})\simeq
K_*(R)\oplus K_*(R/I_{12})\oplus\cdots\oplus K_*(R/I_{n-2\; n-1})$.
Hence
$$K_*(B)\simeq K_*(R)\oplus K_*(R/I_{12})\oplus \cdots\oplus K_*(R/I_{n-2\,n-1})\oplus K_*(R/I_{n-1\, n}).$$
This finishes the proof of Lemma \ref{ktheory1}. $\square$

In particular, as a consequence of Lemma \ref{ktheory1}, we can
strengthen the result in \cite{gkuk} as the following corollary,
here we drop all assumptions on rings and ideals.

\begin{Koro} Let $R$ be an arbitrary ring with identity and $I$ an arbitrary ideal in $R$.
Then, for a ring of the following form $$S=\begin{pmatrix}
  R        &  I^{t_{12}} &  \cdots     & I^{t_{1n}}    \\
  R        &  R          &  \ddots     &  \vdots  \\
  \vdots   & \vdots      & \ddots      &  I^{t_{n-1\; n}}   \\
  R        & R           & \cdots      & R
\end{pmatrix},$$ where $t_{ij}$ are positive integers, we have
$$  K_*(S)\simeq K_*(R)\oplus \bigoplus_{j=2}^nK_*(R/I^{t_{j-1\;
j}}).$$\label{kgroupbygk} \end{Koro}

As a special case of Corollary \ref{kgroupbygk}, we get the
following result of \cite{keating} without the assumption that $_RI$
and $I_R$ are projective.

\begin{Koro} Let $R$ be a ring with identity and $I$ an ideal in $R$.
Suppose that $t_j$ is a positive integers with $t_j\le t_{j+1}$ for
$j=2, \cdots, n-1$. Let $$T = \begin{pmatrix}
  R        &  I^{t_2} &  \cdots     & I^{t_n}    \\
  R        &  R          &  \ddots     &  \vdots  \\
  \vdots   & \vdots      & \ddots      &  I^{t_n}   \\
  R        & R           & \cdots      & R \end{pmatrix}.$$ Then $T$ is a ring and $$ K_*(T)\simeq K_*(R)\oplus \bigoplus_{i=2}^n K_*(R/I^{t_i}).$$
\label{k2}
\end{Koro}

Let us remark that if $I$ is a nilpotent ideal in $R$ with identity
then $K_0(R)\simeq K_0(R/I)$. In general, this is not true for
higher $ K_n$-groups with $n\ge 1$. Thus, for $K_0$, we may replace
the direct summands $K_0(R/I^{t_j})$ by $ K_0(R/I)$ in Corollary
\ref{k2}, and get $ K_0(T)\simeq K_0(R)\oplus (n-1) K_0(R/I).$

Similarly, we have the following result on the groups $K_n$ of the
ring defined in Lemma \ref{de}

\begin{Lem} Let $B$ be the ring
$$B :=
 \begin{pmatrix}
R & I_2 &I_3 & \cdots & I_{n-1} & I_n\\
R & R_2 & I_3 & \cdots & I_{n-1} & I_n\\
R & I_{32}& R_3& \ddots &  \vdots & I_n\\
R & I_{42} & I_{43}& \ddots & I_{n-1} & \vdots\\
\vdots & \vdots &\vdots & \ddots & R_{n-1}& I_n\\
R & I_{n2} & I_{n3} & \cdots & I_{n\, n-1} & R_n
 \end{pmatrix}$$
defined in Lemma \ref{de}. Then $$  K_*(B)\simeq K_*(R)\oplus
\bigoplus_{j=2}^nK_*(R_j/I_j).$$ \label{ktheory2}
\end{Lem}

This result shows that the abelian group $K_n(B)$ of the ring $B$ is
independent of the choice of the ideals $I_{ij}$ in $R$ for all
$n\ge 0$.

As a direct consequence of Lemma \ref{ktheory2}, we have the
following corollary.

\begin{Koro}Let $R$ be a ring with identity, and let $I_j$ be an
ideal of $R$ with $2\le j\le n$ such that $I_j\subseteq I_{j-1}$ for
all $j$. Then, for the rings

$$S:=\begin{pmatrix}
  R & I_2& I_3 & \cdots & I_n\\
  R & R & I_3  & \cdots & I_n\\
  R & I_2 & R  & \ddots & \vdots\\
  \vdots & \vdots & \ddots & \ddots & I_n \\
  R &I_2 & \cdots & I_{n-1} & R
\end{pmatrix}, \qquad T:=\begin{pmatrix}
  R & I_2 & I_3 &  \cdots  & I_n\\
  R & R & I_3 &  \cdots  & I_n\\
  R & R & R &  \ddots  & \vdots\\
  \vdots & \vdots & \vdots &\ddots  &I_n\\
  R & R & R &  \ldots & R
\end{pmatrix},$$
we have $$ K_*(S)\simeq K_*(R)\oplus
\bigoplus_{j=2}^nK_*(R/I_j)\simeq K_*(T).$$ \label{k3}
\end{Koro}

Let us remark that we can also use our method in this section to
calculate some corner rings $eBe$, though, in general, we cannot get
an $\add(_{eBe}eAe)$-split sequence $$0\lra eBe\lra eAe\lra
eAe/eBe\lra 0,$$ with $e$ an idempotent in $B$, from a given
$\add(_BA)$-split sequence $$0\lra B\lra A\lra A/B\lra 0.$$ For
example, suppose that $B$ is the ring defined in Lemma
\ref{ktheory1}. If $e$ is an idempotent in $R$, then, for the corner
ring
$$B_1:=\begin{pmatrix}
eRe & eI_{12}e &eI_{13}e & \cdots & eI_{1\;n}e\\
eRe & eRe & eI_{23}e & \cdots & eI_{2\;n}e \\
\vdots & \vdots &\ddots &\ddots  & \vdots\\
eRe & eRe & \cdots&  eRe & eI_{n-1\;n}e\\
eRe& eRe& \cdots& eRe & eRe
\end{pmatrix}$$ of $B$, we have
$$ K_*(B_1)\simeq K_*(eRe)\oplus \bigoplus_{j=1}^{n-1}K_*(eRe/eI_{j\;j+1}e).$$

Also, we remark that, for any ring $R$, the $R$-duality
$\Hom_R(-,{}_RR)$ is an equivalence between the category $R$-proj
and the category $R^{\opp}$-proj, where $R^{\opp}$ is the opposite
ring of $R$. Thus, for each $n\ge 0$, we have $K_n(R)\simeq
K_n(R^{op})$. From this fact, or from Lemma \ref{lem1}(3) for right
modules, we can see that if $S'$ is a ring of the form
$$S':=\begin{pmatrix}
  R & I_1 & I_1 &  \cdots  & I_1\\
  I_2& R & I_2 &  \cdots  & I_2\\
  \vdots& \ddots &\ddots &  \ddots  & \vdots\\
  I_{n-1} &  \cdots & I_{n-1}&R  &I_{n-1}\\
  R & \cdots  & R&  R & R
\end{pmatrix},$$ where $R$ is a ring with identity and $I_j$ is an
ideal of $R$ for each $1\le j<n$, then
$$ K_*(S')\simeq K_*(R)\oplus \bigoplus_{j=1}^{n-1}K_*(R/I_j).$$
Note that $S'$ is closely related to the ring $S$ in Corollary
\ref{k3} (1).

\medskip
Now, recall that a pullback diagram of rings:
$$(*) \qquad \begin{CD}R@>{f_1}>>R_1\\@V{h_2}VV @VV{h_1}V\\
R_2@>{f_2}>>R_0\end{CD}$$ is called a Milnor square if one of $f_2$
and $h_1$ is surjective.

An example of Milnor squares is the following case: Let $R\subseteq
S$ be an extension of rings with the same identity. If there is an
ideal $J$ of $S$ such that $J\subseteq R$, then there is a canonical
Milnor square
$$\begin{CD}R@>>>S\\@VVV @VVV\\
R/J@>>>S/J.\end{CD}$$

Let $R$ be the product $R_1\times \cdots \times R_n$ of finitely
many rings $R_i$ with $1\le i\le n$. A subdirect product of ring $R$
is a subring $S\subseteq R$ for which each projection $S\ra R_i$
carries $S$ onto $R_i$ for each $i$. In this case we say that the
inclusion $S\subseteq R$ is an inclusion of a subdirect product.

The following lemma is useful and well-known for calculation of
higher $K$-groups of rings.

\begin{Lem} For a given Milnor square $(*)$, the following are true:

$(1)$ {\rm (See  \cite[Theorem 3.3]{milnor})} There is a
Mayer-Vietoris exact sequence:
$$ K_1(R)\lraf{\big((f_1)_*,(h_2)_*\big)}K_{1}(R_1)\oplus
K_{1}(R_2)\lraf{{\small \begin{pmatrix}(h_1)_*\\
-(f_2)_*\end{pmatrix}}}K_{1}(R_0) \lra
K_0(R)\lraf{\big((f_1)_*,(h_2)_*\big)}K_{0}(R_1)\oplus
K_{0}(R_2)\lraf{\begin{pmatrix}(h_1)_*\\
-(f_2)_*\end{pmatrix}}K_{0}(R_0),$$ where $f_*$ denotes the
homomorphism induced by $f$.

$(2)$ {\rm (See  \cite{charney}, \cite[Theorem 5.5]{weibel1})}
Suppose that $(*)$ is a Milnor square of ${\mathbb
Z}/p^m\mathbb{Z}$-algebras, where $p\ge 2$ and $m$ are fixed
positive integers. Let $s$ be a non-zero integer such that $p$
divides $s$. Then there is an exact sequence of $K$-groups, that is,
the Mayer-Vietoris sequence:
$$\begin{array}{rl} \cdots\lra & K_{*+1}(R_1)\otimes_{\mathbb Z}{\mathbb
Z}[\frac{1}{s}]\oplus
K_{*+1}(R_2)\otimes_{\mathbb Z}{\mathbb Z}[\frac{1}{s}]\lraf{\begin{pmatrix}(h_1)_*\\
-(f_2)_*\end{pmatrix}}K_{*+1}(R_0)\otimes_{\mathbb Z}{\mathbb
Z}[\frac{1}{s}] \lra K_*(R)\\
 & \lraf{\big((f_1)_*,(h_2)_*\big)}K_{*}(R_1)\otimes_{\mathbb Z}{\mathbb Z}[\frac{1}{s}]\oplus
K_{*}(R_2)\otimes_{\mathbb Z}{\mathbb Z}[\frac{1}{s}]\lraf{\begin{pmatrix}(h_1)_*\\
-(f_2)_*\end{pmatrix}}K_{*}(R_0)\otimes_{\mathbb Z}{\mathbb
Z}[\frac{1}{s}]\lra\cdots.\end{array}$$

$(3)$ Suppose that $(*)$ is a Milnor square of ${\mathbb
Z}/p^m\mathbb{Z}$-algebras, where $p\ge 2$ and $m$ are fixed
positive integers. Let $s$ be a non-zero integer such that $p$
divides $s$. If the induced homomorphism $(f_2)_*$ in $(2)$ is an
split epimorphism for all $*\in \mathbb N$, then there is an exact
sequence for all $*\in \mathbb N$:
$$\begin{array}{c}0\lra K_*(R)\otimes_{\mathbb Z}{\mathbb Z}[\frac{1}{s}]\lraf{\big((f_1)_*,(h_2)_*\big)}K_{*}(R_1)
\otimes_{\mathbb Z}{\mathbb Z}[\frac{1}{s}]\oplus
K_{*}(R_2)\otimes_{\mathbb Z}{\mathbb Z}[\frac{1}{s}]\lraf{{\small\begin{pmatrix}(h_1)_*\\
-(f_2)_*\end{pmatrix}}}K_{*}(R_0)\otimes_{\mathbb Z}{\mathbb Z}[
\frac{1}{s}]\lra 0.\end{array}$$

In particular, if the induced homomorphism $(f_2)_*$ in $(2)$ is an
isomorphism for all $*\in \mathbb N$, then so is the induced
homomorphism $(f_1)_*$.

 $(4)$ {\rm \cite[Theorem 13.33]{magurn}, \cite{milnor}} If both $h_1$
and $f_1$ are surjective, or if $h_1$ is surjective and $f_1$ is the
inclusion of a subdirect product, then there is an exact sequence
$$\begin{array}{rl} K_2(R)&\lraf{\big((f_1)_*,(h_2)_*\big)}K_{2}(R_1)\oplus
K_{2}(R_2)\lraf{\begin{pmatrix}(h_1)_*\\
-(f_2)_*\end{pmatrix}}K_{2}(R_0)\lra
K_1(R)\lraf{\big((f_1)_*,(h_2)_*\big)}\\ & K_{1}(R_1)\oplus
K_{1}(R_2)\lraf{\begin{pmatrix}(h_1)_*\\
-(f_2)_*\end{pmatrix}}K_{1}(R_0)\lra
K_0(R)\lraf{\big((f_1)_*,(h_2)_*\big)}K_{0}(R_1)\oplus
K_{0}(R_2)\lraf{\begin{pmatrix}(h_1)_*\\
-(f_2)_*\end{pmatrix}}K_{0}(R_0).\end{array}$$ \label{mvs}
\end{Lem}

Remark that there is a dual statement of (3) for $(f_1)_*$ to be a
split monomorphism for each $*\in \mathbb N$.

\medskip
Now we turn to proof of Theorem \ref{thm1}. Observe that the
argument in our proof below is actually a combination of the
previous results with Mayer-Vietoris exact sequences, and works also
for many other cases. Here we prove only Theorem \ref{thm1}(1), and
leave the details of the proof of Theorem \ref{thm1}(2) to the
reader.

\medskip
{\bf Proof of Theorem \ref{thm1}} (1): Let
$$J:=\begin{pmatrix}
  I & I_{12} &   \cdots  & I_{1n}\\
  I & I &  \ddots  & \vdots\\
  \vdots& \ddots  &\ddots    & I_{n-1\, n}\\
  I & \cdots & I  & I
\end{pmatrix}, \qquad
B:=\begin{pmatrix}
  R & I_{12} &   \cdots  & I_{1n}\\
  I& R &   \ddots  & \vdots\\
  \vdots& \ddots &\ddots &   I_{n-1\, n}\\
  I &  \cdots & I  &R\\
  \end{pmatrix}, \qquad
  A:=\begin{pmatrix}
  R & I_{12} &   \cdots  & I_{1n}\\
  R& R &   \ddots  & \vdots\\
  \vdots& \ddots &\ddots &   I_{n-1\,n}\\
  R &  \cdots & R  &R\\
  \end{pmatrix}.$$
By the assumptions in Theorem \ref{thm1}(1), we can verify that $A$
and $B$ are rings and that $J$ is an ideal of $A$. Thus $J$ is also
an ideal of $B$. Note that $B$ is a subalgebra of $A$. Let $f$ be
the inclusion of $B$ into $A$. If we define

$$B':=B/J=\begin{pmatrix}
  R/I & 0 &   \cdots  & 0\\
  0 & R/I &   \ddots  & \vdots\\
  \vdots& \ddots &\ddots &   0\\
  0 &  \cdots & 0  &R/I\\
  \end{pmatrix}, \qquad
  A':=A/J=\begin{pmatrix}
  R/I & 0 &   \cdots  & 0\\
  R/I& R/I&   \ddots  & \vdots\\
  \vdots& \ddots &\ddots &   0\\
  R/I &  \cdots & R/I  &R/I\\
  \end{pmatrix},$$  then we have a Milnor square

$$\begin{CD} B @>f>> A\\ @V{g'}VV @VV{g}V\\
B'@>f'>> A',\end{CD}$$ where $g$ and $g'$ are the canonical
surjective maps, and where $f'$ is the injective map induced from
$f$. Since the map $f'_*: K_*(B')\ra K_*(A')$ is an isomorphism for
$*=0,1, 2,\cdots,$ it follows that $f'_*\otimes {\mathbb
Z}[\small{\frac{1}{s}}]: K_*(B')\otimes {\mathbb
Z}[\small{\frac{1}{s}}]\ra K_*(A')\otimes {\mathbb
Z}[\small{\frac{1}{s}}]$ is an isomorphism. Thus we see from Lemma
\ref{mvs}(3) that $f_*\otimes {\mathbb Z}[\small{\frac{1}{s}}]$ is
also an isomorphism. It then follows from Corollary
\ref{ktheory1}(2) that
$$\begin{array}{c}K_*(B)\otimes_{\mathbb Z}{\mathbb Z}[\frac{1}{s}]\simeq K_*(A)\otimes_{\mathbb Z}{\mathbb Z}[\frac{1}{s}]\simeq
K_*(R)\otimes_{\mathbb Z}{\mathbb Z}[\frac{1}{s}]\oplus
\bigoplus_{j=1}^{n-1}K_*(R/I_{j\, j+1})\otimes_{\mathbb Z}{\mathbb
Z}[\frac{1}{s}].\end{array}$$ This finishes the proof of Theorem
\ref{thm1}(1).

If we define $$J:=\begin{pmatrix}
  I & I_{2} &   I_3 &\cdots  & I_{n}\\
  I & I_2 &  I_3 &\cdots  & \vdots\\
  I & I_{32} & I_3 &\ddots &\vdots\\
  \vdots& \vdots &\ddots  &\ddots    & I_{n}\\
  I & I_{n2}& \cdots & I_{n\, n-1}  & I_n\\
\end{pmatrix},$$ then the proof of Theorem \ref{thm1}(2) can be
carried out similarly since we have Lemma \ref{ktheory2}. $\square$

\medskip
Now we mention the following corollary of the proof of Theorem
\ref{thm1}.

\begin{Koro} Suppose that $p\ge 2$ and  $m$ are positive integers.
Let $R$ be a ${\mathbb Z}/p^m\mathbb{Z}$-algebra with identity, and
let $I$ and $J$ be two arbitrary ideals of $R$. Define $$S:=
{\begin{pmatrix}
  R & I &   \cdots  & I\\
  J& R &   \ddots  & \vdots\\
  \vdots& \ddots &\ddots &   I\\
  J &  \cdots & J  &R\\
\end{pmatrix}.}_{n\times n}$$ Then $S$ is a ring, and if, in addition, $I^2\subseteq J$ (for example, $I^2=0$, or $I\subseteq J$),
we have
$$\begin{array}{c} K_*(S)\otimes_{\mathbb Z}{\mathbb Z}[\frac{1}{s}]\simeq  K_*(R)\otimes_{\mathbb Z}{\mathbb Z}[\frac{1}{s}]
\oplus (n-1)K_*(R/I)\otimes_{\mathbb Z}{\mathbb
Z}[\frac{1}{s}]\end{array}$$ for all non-zero integer $s$ such that
$s$ is divisible by $p$.\label{k6}
\end{Koro}

{\it Proof.} We define
$$A:= \begin{pmatrix}R & I &   \cdots  & I\\
I+ J& R &   \ddots  & \vdots\\
\vdots& \ddots &\ddots &   I\\
I+J &  \cdots & I+ J  &R\\
\end{pmatrix}_{n\times n},\qquad
J':= \begin{pmatrix}I & I &   \cdots  & I\\
J& I &   \ddots  & \vdots\\
\vdots& \ddots &\ddots &   I\\
J &  \cdots & J  &I\\
\end{pmatrix}_{n\times n}.$$ Then one can verify that $A$ is a ring and $J'\subseteq S$ is
an ideal in $A$. Note that $S$ is a subring of $A$. Now, let $B:=S$,
$B':=B/J'$ and $A':=A/J'$. Then we may use the same argument as in
the proof of Theorem \ref{thm1} to show $K_*(B)\otimes_{\mathbb
Z}{\mathbb Z}[\small{\frac{1}{p}}]\simeq K_*(A)\otimes_{\mathbb
Z}{\mathbb Z}[\small{\frac{1}{p}}]$. But for the latter, we have
$K_*(A)\otimes_{\mathbb Z}{\mathbb Z}[\small{\frac{1}{p}}]\simeq
K_*(R)\otimes_{\mathbb Z}{\mathbb Z}[\small{\frac{1}{p}}]\oplus
(n-1)K_*(R/I)\otimes_{\mathbb Z}{\mathbb Z}[\small{\frac{1}{p}}]$ by
Theorem \ref{thm1}(1). Thus Corollary \ref{k6} follows. $\square$

\medskip
Let us illustrate how the argument in the above proof of Theorem
\ref{thm1}(1) can be applied to other cases.

Again, suppose that $p\ge 2$ and $m$ are positive integers. Let $R$
be a ${\mathbb Z}/p^m\mathbb{Z}$-algebra with identity and $I$ an
arbitrary ideal of $R$. For each finite partially ordered set $P$,
we associate a ring $B:=B(R,I,P)$ which is a subring of the matrix
ring over $R$ with indexing set $P$, it is defined as follows: Let
$B=(B_{ij})_{i,j\in P}$ with $ B_{ij}=R$ if $i\ge j$, and $B_{ij}=I$
otherwise. We may assume that $P=\{a_1, \cdots, a_n\}$ such that
$a_i\le a_j$ implies $i\le j$. Under this assumption we see that
$J':=M_n(I)$ is an ideal of $B$, which is also an ideal of
$$ A:={\begin{pmatrix}
  R & I &   \cdots  & I\\
  R& R &   \ddots  & \vdots\\
  \vdots& \ddots &\ddots &   I\\
  R &  \cdots & R  &R\\
  \end{pmatrix}.}_{n\times n}$$
Note that $B$ is a subring of $A$. Let $B':=B/J'$ and $A':=A/J'$. We
define $C$ to be the diagonal matrix ring with the principal
diagonal entries $R/I$. Then $C$ is a subring of both $B'$ and $A'$.
Using this ring $C$, we can see that the inclusion $f'$ of $B'$ into
$A'$ induces an isomorphism $(f')_*$ for all $*\in \mathbb N$. Then
we may use the same argument as the above to show that, for any $s$
divisible by $p$,
$$ \begin{array}{c} K_*(B)\otimes_{\mathbb
Z}{\mathbb Z}[\small{\frac{1}{s}}]\simeq K_*(A)\otimes_{\mathbb
Z}{\mathbb Z}[\small{\frac{1}{s}}]\simeq K_*(R)\otimes_{\mathbb
Z}{\mathbb Z}[\small{\frac{1}{s}}]\oplus
(n-1)K_*(R/I)\otimes_{\mathbb Z}{\mathbb
Z}[\small{\frac{1}{s}}].\end{array}$$

\medskip
We end this section by a couple of remarks concerning Theorem
\ref{thm1}.

(1) In Theorem \ref{thm1}, if $R$ is a $\mathbb{Z}_p$-algebra
instead of a $\mathbb{Z}/p^m\mathbb{Z}$-algebra, and if $R/I$,
$R/I_i$ and $R/I_{ij}$ are finite rings for all $i, j$, then Theorem
\ref{thm1} still holds true. Indeed, in this case we can use
Charney's excision at the end of the paper \cite{charney} since
$I\otimes_{\mathbb{Z}}\mathbb{Z}[\frac{1}{s}]$ has a unit. This is
due to $\Tor_1^{\mathbb{Z}}(-,\mathbb{Z}[\frac{1}{s}])=0$ and to the
fact that the quotient rings $R/I$, $R/I_i$ and $R/I_{ij}$ are
$\mathbb{Z}/p^m\mathbb{Z}$-algebras for some $m>0$. Indeed, we have
an exact sequence $$\Tor^{\mathbb{Z}}_1(R/I,
\mathbb{Z}[\frac{1}{s}])\lra
I\otimes_{\mathbb{Z}}\mathbb{Z}[\frac{1}{s}]\lra
R\otimes_{\mathbb{Z}}\mathbb{Z}[\frac{1}{s}]\lra
(R/I)\otimes_{\mathbb{Z}}\mathbb{Z}[\frac{1}{s}].$$ Clearly, the
first and last terms vanish, this implies
$I\otimes_{\mathbb{Z}}\mathbb{Z}[\frac{1}{s}]\simeq
R\otimes_{\mathbb{Z}}\mathbb{Z}[\frac{1}{s}]$. So, the condition of
Charney's result in \cite{charney} is satisfied. I thank X. J. Guo
for explanation of this fact.

(2) A crucial fact of our proofs of the main results is: Given an
extension $B\subseteq A$ of rings with the same identity such that
$_BA$ is finitely generated and projective, we have $K_*(B)\simeq
K_*(\End_B(A\oplus \, A/B))$ for all $*\in \mathbb{N}$. Moreover, we
may also compare the algebraic $K$-theory of $B$ with that of $A$.
For this purpose, we define $\Omega$ to be the kernel of the
multiplication map $A\otimes_BA\ra A$, it follows from the
Additivity Theorem (see \cite[Corollary 1, Section 3]{Quillen} that
the exact sequence of the exact functors
$$ 0\lra \Omega\otimes_A-\lra A\otimes_B-\lra id \lra 0$$
on the category of finitely generated projective $A$-modules gives
rise to three homomorphisms of abelian groups: $r_*: K_*(A)\ra
K_*(B)$, $t_*: K_*(B)\ra K_*(A)$ and $\omega_*: K_*(A)\ra K_*(A)$
such that $r_*t_*=1_{K_*(A)} +\omega_*$. If, in addition, the
$n$-fold tensor product of $\Omega$ over $A$ vanishes for some
natural number $n$, that is, $\Omega^{{\otimes_A}n}=0$ (for example,
$\Omega = 0$ in case the inclusion $B\subseteq A$ is an injective
ring epimorphism), then the map $t_*$ is split surjective, and
$K_*(A)$ is a direct summand of $K_*(B)$. In general, neither $t_*$
nor $r_*$ is an isomorphism.

\section{Lower $K$-theory for matrix subrings \label{lowerK}}

In this section we consider the algebraic $K$-groups $K_0$ and $K_1$
for matrix subrings. Our results in this section are not covered by
the main results in the previous sections.

We first consider the group $K_0$. In this case, we have the
following result in which we do not assume that the rings considered
are $\mathbb{Z}/p^m\mathbb{Z}$-algebras or
$\mathbb{Z}[\frac{1}{p}]$-algebras.

\begin{Prop} Let $R$ be an arbitrary ring with identity, and let $I, J$ and $I_{ij}$ be ideals in $R$.

$(1)$  For the rings $S$ and $T$ defined in Theorem \ref{thm1}, we
have
$$ K_0(S)\simeq K_0(R)
\oplus \bigoplus_{j=2}^nK_0(R/I_{j-1\; j}), \quad K_0(T)\simeq
K_0(R) \oplus \bigoplus_{j=2}^nK_0(R_j/I_j).$$

$(2)$ For the ring $S$ defined in Corollary \ref{k6} with
$I^2\subseteq J$, we have
$$K_0(S)\simeq K_0(R)\oplus (n-1)K_0(R/I).$$
\label{k7}
\end{Prop}

The proof of this proposition is actually a combination of Corollary
\ref{k3} and Lemma \ref{mvs}(1) and (3), and we leave the details of
the proof to the interested reader.

Here arises an open question: We do not know, at moment, if
Proposition \ref{k7} is true for higher algebraic $K$-groups $K_n$
with $n\ge 1.$ But, for $K_1$, we do have some partial answers.
Before stating our result, we first prove the following lemma.

\begin{Lem} Let $B\subseteq A$ be an extension of rings with the same identity.
Suppose that $I$ is an idempotent ideal of $A$ contained in $B$. If
the inclusion $B\subseteq A$ induces an isomorphism $\gamma_i:
K_i(B/I)\ra K_i(A/I)$ for $i=1,2,$ then $K_1(B)\simeq K_1(A)$.
\label{k8}
\end{Lem}

{\it Proof.} Let $K_i(B,I)$ denote the $i$-th relative $K$-group of
the canonical surjective map $B\ra B/I$. Then we may form the
following commutative diagram of ableian groups with exact rows:
$$\begin{CD} K_2(B/I)@>>>K_1(B,I)@>>>K_1(B)@>>>K_1(B/I)@>>>K_0(B,I)\\
@VV{\gamma_2}V @VV{\gamma}V
@VV{\beta}V@VV{\gamma_1}V@VV{\simeq}V\\
K_2(A/I)@>>>K_1(A,I)@>>>K_1(A)@>>>K_1(A/I)@>>>K_0(A,I).
\end{CD}$$
Here we use the fact that $K_0(B,I)$ is always independent of $B$.
Now, by the Five Lemma in homological algebra, we know that the map
$\beta$ is isomorphic if $\gamma$ is isomorphic. However, this
follows from a result of Vaserstein (see \cite[Chapter III, Section
2, Remark 2.2.1]{weibel}), which states that if $J$ is an ideal in a
ring $R$ with identity, then $K_1(R,J)$ is independent of $R$ if and
only if $J^2=J$.  Thus $\gamma$ is an isomorphism. $\square$

\medskip
We should notice that, in general, $K_n(R,I)$ depends on $R$ for
$n\ge 1$. This is why the conclusions in Theorem \ref{thm1} are
localized.

So, with Lemma \ref{k8} in hand, we can prove the following
proposition for $K_1$.

\begin{Prop} Let $R$ be a ring with identity, and let $I\subseteq J$ be ideals in $R$. If $I$ is idempotent, then, for the ring
$$B:=\begin{pmatrix}
  R & I &   \cdots  & I\\
  J& R &   \ddots  & \vdots\\
  \vdots& \ddots &\ddots &   I\\
  J&  \cdots & J  &R\\
\end{pmatrix}_{n\times n}, $$ we have
$$K_1(B)\simeq K_1(R)\oplus (n-1)K_1(R/I).$$
\label{k9}
\end{Prop}

{\it Proof.} Clearly,  $B$ is a subring of the ring
$$A:=\begin{pmatrix}
  R& I &   \cdots  & I\\
  R& R &   \ddots  & \vdots\\
  \vdots& \ddots &\ddots &   I\\
  R &  \cdots & R  &R\\
\end{pmatrix}_{n\times n},$$ and  $J':=M_n(I)$, the $n\times n$ matrices over $I$, is an idempotent ideal of
$A$ and $B$, respectively. We know that $K_*(B/J')$ and $K_*(A/J')$
are isomorphic for all $*\in \mathbb{N}$. Hence Proposition \ref{k9}
follows from Lemma \ref{k8} and Corollary \ref{ktheory1}
immediately. $\square$

Finally, we mention another type of matrix rings: Let $R$ and $S$
are rings with identity and $_RM_S$ and $_SN_R$ are bimodules. We
define a ring $${A:=\begin{pmatrix}
  R& M \\
  N& S \\
  \end{pmatrix}},\quad{\begin{pmatrix}
  r& m \\
  n& s \\
  \end{pmatrix}\cdot\begin{pmatrix}
  r'& m' \\
  n'& s' \\
  \end{pmatrix}=\begin{pmatrix}
  rr'& rm'+ms' \\
  nr'+sn'& ss' \\
  \end{pmatrix}}$$
for $r,r'\in R, s,s'\in S, m,m'\in M$ and $n,n'\in N$. Note that
${M':=\begin{pmatrix}
  0& M \\
  0& 0 \\
  \end{pmatrix}}$ and $N':={\begin{pmatrix}
  0& 0 \\
  N& 0 \\
  \end{pmatrix}}$
are two ideals in $A$. Thus one has a Milnor diagram $$
\begin{CD}A@>>>A/M'\\@VVV @VVV \\ A/N' @>>> A/(M'+N').\end{CD}$$
By Lemma \ref{mvs}(4), we can show that $K_i(A)\simeq K_i(R)\oplus
K_i(S)$ for $i=0,1.$ This result can be used to reduce the
calculation of lower $K$-groups of finite-dimensional algebras with
radical-square-zero to local algebras.

\section{Higher mod-$p$ $K$-theory\label{modp}}

In this section, we shall point out that our main result, Theorem
\ref{thm1}, holds true for the mod-$p$ $K$-theory $K_*(-,
\mathbb{Z}/p\mathbb{Z})$ under the assumption that algebras
considered are $\mathbb{Z}[\frac{1}{p}]$-algebras and $p\not\equiv 2
(\mbox{mod}\, 4)$, where $p\ge 2$ is any positive integer.

Let $R$ be a ring with identity.  In \cite{browder}, Browder
developed $K$-theory with coefficients $\mathbb{Z}/p\mathbb{Z}$.
This is the so-called mod-$p$ $K$-theory
$K_*(R,\mathbb{Z}/p\mathbb{Z})$ for $*\in {\mathbb Z}$. Note that
$K_0(R,\mathbb{Z}/p\mathbb{Z})=K_0(R)\otimes_{\mathbb{Z}}\mathbb{Z}/p\mathbb{Z}$,
and $K_i(R, \mathbb{Z}/p\mathbb{Z})=0$ if $i<0$ (see \cite [p.
45]{browder}). Later, Weibel observed in \cite{weibel2} that
excision holds and that Mayer-Vietoris sequences exist if the rings
involved are $\mathbb{Z}[\frac{1}{p}]$-algebras. The mod-$p$
$K$-theory is closely related to the usual $K$-theory in the
following manner.

\begin{Lem} {\rm Universal Coefficient Theorem (see \cite{browder} and \cite{weibel2}):}

Let  $R$ be a $\mathbb{Z}[\frac{1}{p}]$-algebra with identity. For
all $*\in\mathbb{N}$, there is a short exact sequence of abelian
groups
$$ 0\lra K_*(R)\otimes_{\mathbb Z}\mathbb{Z}/p\mathbb{Z}\lra
K_*(R,\mathbb{Z}/p\mathbb{Z})\lra
\Tor_1^{\mathbb{Z}}(K_{*-1}(R),\mathbb{Z}/p\mathbb{Z})\lra 0.$$ If
$p\not\equiv 2 \, \emph{(mod 4)}$, then this sequence splits (not
naturally), so that $K_*(R,\mathbb{Z}/p\mathbb{Z})$ is a
$\mathbb{Z}/p\mathbb{Z}$-module. If $p\equiv 2 \; \emph{(mod 4)}$,
then $K_*(R,\mathbb{Z}/p\mathbb{Z})$ is a
$\mathbb{Z}/2p\mathbb{Z}$-module. \label{uct}
\end{Lem}

Thus, if $p\not\equiv 2 \, \mbox{(mod 4)}$, we see that
$K_*(R,\mathbb{Z}/p\mathbb{Z})$ is completely determined by the
usual $K$-groups $K_*(R)$.

Another result which we need is a Mayer-Vietoris sequence for
mod-$p$ $K$-groups.

\begin{Lem} {\rm \cite[Corollary 1.3]{weibel2}} For a Milnor square
$(*)$ of $\,\mathbb{Z}[\frac{1}{p}]$-algebras, there is a long exact
sequence of abelian groups for all integers $*:$

$$ \begin{array}{rl} \cdots \lra  & K_{*+1}(R_1, \mathbb{Z}/p\mathbb{Z})\oplus  K_{*+1}(R_2, \mathbb{Z}/p\mathbb{Z})\lra
K_{*+1}(R_0,\mathbb{Z}/p\mathbb{Z})\lra \\  & \\ & K_*(R,
\mathbb{Z}/p\mathbb{Z}) \lra K_{*}(R_1,
\mathbb{Z}/p\mathbb{Z})\oplus  K_{*}(R_2,
\mathbb{Z}/p\mathbb{Z})\lra K_*(R_0, \mathbb{Z}/p\mathbb{Z})\lra
\cdots.\end{array}$$
\end{Lem}

Now it follows from the above two lemmas and Theorem \ref{thm2} that
Theorem \ref{thm1} holds true for the mod-$p$ $K$-groups $K_*(R,
\mathbb{Z}/p\mathbb{Z})$ if $R$ is a
$\mathbb{Z}[\frac{1}{p}]$-algebra and if $p\not\equiv 2 \;
\mbox{(mod\, 4)}$, since the argument there in the proof of Theorem
\ref{thm1} works in our new situation.

If $p\equiv 2 \mbox{ (mod 4)}$, we do not know whether
$K_*(R,\mathbb{Z}/p\mathbb{Z})$ can be fully controlled by the first
and last terms in Lemma \ref{uct}. In general, extensions of fixed
abelian groups may not be isomorphic, for instance, the cyclic group
of order $4$  and the Klein group $(\mathbb{Z}/2\mathbb{Z})\times
(\mathbb{Z}/2\mathbb{Z})$ both are extensions of the cyclic group of
order $2$ by itself, but they are not isomorphic.

\section{Examples: GV-ideals\label{example}}

In this section we shall give some examples related to our results.
The first one is constructed from a $\mathcal D$-split sequence
which is induced by a surjective ring homomorphism.

Let $B$ be a ring with identity and $J$ an ideal of $B$. We define
$A=B/J$. Then we have an exact sequence in $B$-Mod:
$$0\ra J\ra B\lraf{\pi}
A\ra 0,$$ where $\pi$ is the canonical surjection.

For this sequence to be an $\add(_BB)$-split sequence, we have to
assume $\Ext_B^1(A,B)=0$. This happens often in commutative algebra.
For example, if $B$ is a commutative noetherian ring, and $J$ is an
ideal of $B$ such that $J$ contains a regular sequence on $B$ of
length $2$, then $\Ext^i_B(A,B)=0$ for $i= 0,1$ (see \cite[p.
101]{kaplansky}). Another example is the so-called GV-ideals in
integral domains. Here we will state the following general
definition of GV-ideals.

Let $R$ be an arbitrary ring with identity. Recall that an ideal $I$
of $R$ is called a GV-ideal (after the names Glaz and Vasconcelos,
see \cite{ywzc, GV}) if the induced map $\mu_I: R\lra \Hom_R(I,R)$,
given by $r\mapsto (x\mapsto xr)$ for $x\in I$, is an isomorphism of
$R$-bimodule. This is equivalent to $\Ext^i_R(R/I,R)=0$ for $i=0,1.$
Thus $R$ is a GV-ideal of $R$. Note that $p\mathbb{Z}$ is not a
GV-ideal of $\mathbb Z$ for any $p\in\mathbb{Z}$ with $|p|\ne 1$,
even though we have ${\mathbb Z}\simeq \Hom_{\mathbb
Z}(p\mathbb{Z},\mathbb{Z})$. We remark that the above definition of
GV-ideals is more general than that in commutative rings where it is
required that $_RI$ is finitely generated (see \cite{ywzc}).

Let $GV(R)$ be the set of all GV-ideals of $R$. For ideals $I$ and
$J$ of $R$, we denote by $(I:J):=\{x\in R\mid Ix\subseteq J\}.$
(This notation is different from what was usually used in ring
theory, but soon we will see its convenience when elements compose).
Clearly, $(I:R)=R, (R:I)=I$, and $(I:J)$ is an ideal of $R$.

The following lemma shows some properties of GV-ideals, which are of
interest for our proofs.

\begin{Lem} Let $B$ be a ring with identity, and let $J$ be a GV-ideal in
$B$. Then

$(1)$ the sequence $0\ra J\ra B\lraf{\pi} A\ra 0$ is an
$\add(_BB)$-split sequence in $B\Modcat$. Thus $\End_B(B\oplus J)$
is derived-equivalent to $\begin{pmatrix}
B & B/J  \\
0 & B/J \end{pmatrix}.$

$(2)$ $\End_B(J)\simeq B$ $($as rings and as $B$-bimodules$)$.

$(3)$ If $I$ is an ideal in $B$, then $_B\Hom_B(J,I)_B\simeq (J:I)$
as $B$-bimodules. In particular, if $J\subseteq I$, then
$\Hom_B(J,I)\simeq B$.

$(4)$ If $x\in B$ such that $Jx=0$, then $x=0$.

$(5)$ If $I$ is an ideal in $B$ with $J\subseteq I$, then $I\in
GV(B)$.

$(6)$ If $I,J\in GV(B)$, then $IJ\in GV(B)$.\label{gv1}
\end{Lem}

{\it Proof.} (1) is clear by \cite[Theorem 1.1]{hx2} and the
definition of GV-ideals.

(2) By the definition of GV-ideals, the induced map $\mu_J: B\lra
\Hom_B(_BJ,B)$ is an isomorphism, this means that every homomorphism
$f$ from $_BJ$ to $_BB$ is given by the right multiplication of an
element in $B$. Since $J$ is an ideal in $R$, $f$ is in fact an
endomorphism of the module $_BJ$. Conversely, if $f\in \End_B(J)$,
then $f$ is a restriction of a right multiplication of an element of
$B$. Hence $\End_B(J)\simeq B$.

(3) We define a map $\varphi: \Hom_B(J,I)\ra (J:I)$ as follows: For
$f\in \Hom_B(J,I)$, there is a unique element $b\in B$ such that the
composition of $f$ with the inclusion $\lambda: I\ra B$ is the right
multiplication map $\cdot b$ since $J\in GV(B)$. This means that
$f\lambda = \cdot b$ and $b\in (J:I)$. So, we define
$f\stackrel{\varphi}{\mapsto} b$. As $(J:I)$ is an ideal of $B$, it
has a canonical bimodule structure. Now one can check that $\varphi$
is an isomorphism of $B$-bimodules.

(4) This is a trivial consequence of the induced isomorphism $\mu_J:
B\simeq \Hom_B(J,B)$.

(5)-(6) These statements were already proved in detail in
\cite{ywzc} for commutative rings, the ideas of their proofs are as
follows: It follows from (4) that $\Hom_B(I/J,B)=0$. Further, by the
isomorphism $\mu_J$ and the fact $\mu_J=\mu_I i_*$ where $i_*:
\Hom_B(I,B)\ra \Hom_B(J,B)$ is induced from the inclusion $i: J\ra
I$, one can check that $\mu_I$ is an isomorphism of $B$-bimodules.
This proves (5).

Let $I,J\in GV(B)$. It follows from (4) that $\mu_{IJ}: B\ra
\Hom_B(IJ,B)$ is injective. We show that it is also surjective. In
fact, since the composition of the maps $B\ra\Hom_B(J,B)\ra
\Hom_B(J,\Hom_B(I,B))\ra \Hom_B(I\otimes_BJ,R)$ is an isomorphism of
$B$-$B$-bimodules, which is the composition of $\mu_{IJ}$ with the
injective map $m_*: \Hom_B(IJ,B)\ra \Hom_B(I\otimes_BJ,B)$ induced
from the surjective multiplication map $I\otimes_BJ\ra IJ$, we see
that $m_*$ is surjective, thus it is an isomorphism of
$B$-$B$-bimodules. This implies that $\mu_{IJ}$ is surjective, and
therefore (6) holds. $\square$

From Lemma \ref{gv1}, we have the following

\begin{Prop} Let $B$ be a ring with identity. Suppose
that $I_n\subseteq I_{n-1}\subseteq \cdots \subseteq I_2\subseteq
I_1$ is a chain of ideals in $B$. If $I_n$ is a GV-ideal in $B$,
then

$(1)$ $\End_B(I_1\oplus \cdots \oplus I_n)$ is isomorphic to
$$C:=\begin{pmatrix}
  B       & (I_1:I_2)      &  (I_1:I_3)        & \cdots  & (I_1:I_n)\\
  B       & B        &  (I_2:I_3)  & \cdots  & (I_2:I_n)\\
  \vdots  & \vdots   &\ddots       & \ddots  &  \vdots\\
  B  & B   &\cdots       & B  & (I_{n-1}:I_n)\\
  B       &  B  & \cdots       & B      & B\\
  \end{pmatrix}. $$

$(2)$ $K_*\big(\End_B(\bigoplus_{j=1}^{n}I_j)\big)\simeq
K_*(B)\oplus \bigoplus_{j=1}^{n-1}K_*\big(B/(I_j:I_{j+1})\big)$ for
all $*\in \mathbb{N}$. \label{gv2}
\end{Prop}

{\it Proof.} (1) Note that $\End_B(I_1\oplus I_2\oplus \cdots \oplus
I_n)$ is the matrix ring with the entries $\Hom_B(I_i,I_j)$ for
$1\le i,j\le n$. Since $I_n$ is a GV-ideal in $B$, every ideal $I_j$
in the chain is a GV-ideal of $B$ by Lemma \ref{gv1}(5). Now (1)
follows from Lemma \ref{gv1} immediately.

(2) This is a direct consequence of (1) and Lemma \ref{ktheory1}.
$\square$

\medskip
As a consequence of Proposition \ref{gv2} and Lemma \ref{gv1}(6), we
have the following corollary.

\begin{Koro} If $I$ is a GV-ideal in a ring $B$ with
identity, then, for any positive integer $n$,
$$ K_*\big(\End_B(\bigoplus_{j=1}^nI^j)\big)\simeq K_*(B)\oplus \bigoplus_{j=1}^{n-1} K_*\big(B/(I^j:I^{j+1})\big).$$
\end{Koro}

If we take $I_1=B$, we have the following corollary.

\begin{Koro} Let $B$ be a ring with identity. Suppose
that $I_n\subseteq I_{n-1}\subseteq \cdots \subseteq I_2\subseteq
I_1=B$ is a chain of $GV$-ideals in $B$. Then
$$ K_*\big(\End_B(B\oplus \bigoplus_{j=2}^nB/I_j)\big)\simeq K_*(B)\oplus K_*(B/I_2)\oplus \bigoplus_{j=2}^{n-1} K_*\big(B/(I_j:I_{j+1})\big).$$
\end{Koro}

{\it Proof.} For each $j$, we have an $\add(_BB)$-split sequence by
Lemma \ref{gv1}(1):
$$ 0\lra I_j\lra {}_BB\lra B/I_j\lra 0.$$
This yields another $\add(_BB)$-split sequence
$$  0\lra \bigoplus_{j=1}^nI_j\lra \bigoplus_{j=1}^n{}_BB\lra \bigoplus_{j=1}^nB/I_j\lra 0.$$
Hence $\End_B(_BB\oplus \bigoplus_{j=2}^nI_j)$ and $\End_B(_BB\oplus
\bigoplus_{j=2}^nB/I_j)$ are derived-equivalent by \cite[Theorem
1.1]{hx2}, and have the isomorphic algebraic $K$-groups $K_*$. By
Lemma \ref{gv2}, we see that $K_*\big(\End_B(_BB\oplus
\bigoplus_{j=1}^nB/I_j)\big)\simeq K_*(B)\oplus
\bigoplus_{j=1}^{n-1}K_*\big(B/(I_j:I_{j+1})\big)$ for all $*\in
\mathbb{N}$. $\square$

\medskip
As a concrete example, we consider the polynomial ring
$B$:=${\mathbb Z}[x]$ over $\mathbb{Z}$ in one variable $x$ and its
ideal $J:=(p,x)$ with $p$ a prime number in $\mathbb{Z}$. It is
known that $J$ is a GV-ideal in $B$. Thus, for the ring
$R:=\End_{{\mathbb Z}[x]}({\mathbb Z}[x]\oplus J)$, by Proposition
\ref{gv2}, we have
$$ K_*(R)\simeq K_*({\mathbb Z}[x])\oplus K_*({\mathbb
Z}/p\mathbb{Z}).$$ Since $\mathbb{Z}$ is a left noetherian ring of
global dimension one, the Fundamental Theorem in algebraic
$K$-theory says that the above isomorphism can be rewritten as
$$ K_*(R)\simeq
K_*({\mathbb Z})\oplus K_*({\mathbb Z}/p\mathbb{Z}). $$ By
\cite{qff}, we get
$$K_0(R)\simeq \mathbb{Z}\oplus \mathbb{Z},\qquad K_1(R)\simeq \mathbb{Z}/2\mathbb{Z}\oplus (\mathbb{Z}/p\mathbb{Z})^{\times},$$
$$K_{2m}(R)=K_{2m}(\mathbb{Z})\quad \mbox{for}\;\; m\ge 1, \qquad
K_{2m-1}(R)\simeq K_{2m-1}(\mathbb{Z})\oplus
\mathbb{Z}/(p^{m}-1)\mathbb{Z}\quad \mbox{for}\;\; m\ge 2,$$ where
$(\mathbb{Z}/p\mathbb{Z})^{\times}$ denotes the set of all non-zero
elements of $\mathbb{Z}/p\mathbb{Z}$. Note that $J$ is not a
projective $\mathbb{Z}[x]$-module. In fact, we have a non-split
exact sequence
$$ 0\lra \mathbb{Z}[x]\lraf{\lambda} \mathbb{Z}[x]\oplus \mathbb{Z}[x]\lraf{\pi} J\lra
0,$$ where $\lambda$ sends $f(x)$ to $(xf(x), -pf(x))$, and $\pi$
sends $(f(x),g(x))$ to $pf(x)+xg(x)$ for all $f(x),g(x)\in
\mathbb{Z}[x]$. So, the result in \cite{keating} cannot be applied
to $R$. However, the one in this note is applicable.

\medskip
Finally, we mention the radical-full extensions in \cite{xi2}.
Recall that an  extension $B\subseteq A$ of rings with the same
identity is said to be left radical-full if $\rad(B)$ is a left
ideal of $A$ and $\rad(A)=\rad(B)A$, where $\rad(A)$ stands for the
Jacobson radical of $A$. So, given a left radical-full extension
$B\subseteq A$ of rings, we may form the ring $C:=\begin{pmatrix}
  A & \rad(B) \\
  A & B \end{pmatrix}. $ It follows from our results in this note that
  $K_n(C)\simeq K_n(A)\oplus K_n(B/\rad(B))$ for all $n\ge 0$ since for any ring extension $S\subseteq R$ and any ideal $I$ in
$S$, if $I$ is a left ideal in $R$ then the rings $\begin{pmatrix}
R & I  \\
R & S \end{pmatrix}$ and $\begin{pmatrix}
  S/I & 0 \\
  R/I &  R \end{pmatrix}$ are derived-equivalent by Lemma \ref{de}.

Related to the last example, we have the following open question:

{\bf Question:} Suppose that $I$ and $J$ are two arbitrary ideals in
a ring $R$ with identity. For the ring $S:=\begin{pmatrix}
R & I  \\
J& R \end{pmatrix}$, can one give a  formula for $K_n(S)$ similar to
the one in Theorem \ref{thm2} for $n\in \mathbb{N}$? (See also the
question mentioned in Section \ref{lowerK}).

\bigskip
{\bf Acknowledgements.} The author would like to thank Xuejun Guo
from Nanjing University for many email correspondences and helpful
explanations on results in $K$-theory, and Fanggui Wang from Sichuan
Normal University for discussion on GV-ideals.

\bigskip
June 23, 2011

\medskip
{\footnotesize
\begin{thebibliography}{99}

\bibitem{bk}{{\sc A. J. Berrick} and {\sc M. E. Keating}, The $K$-theory of triangular matrix rings,  \emph{Contemp. Math.} \textbf{55} (1986), Part I,
69-74. }

\bibitem{BB} {{\sc S. Brenner} and {\sc M. R. Butler}, Generalizations of the Bernstein-Gelfand-Ponomarev reflection functors, In:
\emph{Representation theory} II. (Eds: V. Dlab and P. Gabriel),
Lecture Notes in Math.832, 103-169. Springer-Verlag, Berlin, 1980.}

\bibitem{browder}{{\sc W. Browder,} Algebraic $K$-theory with coefficients $\mathbb{Z}/p$, In: \emph{Geometric applications of homotopy theory} (Evanston, 1977), Lecture Notes in Mathematics 657,
40-84. Springer-Verlag, Berlin, 1978.}

\bibitem{charney}{{\sc R. M. Charney}, A note on excision in $K$-theory, In: \emph{Algebraic
K-theory, number theory, geometry and analysis} (Bielefeld, 1982),
Lecture Notes in Math. 1046, 47-54. Springer-Verlag, Berlin, 1984.}

\bibitem{chyp} {{\sc Y. P. Chen,} Constructions of derived equivalences, Ph. D. Dissertation, 2011.}

\bibitem{DS}{{\sc D. Dugger} and {\sc B. Shipley}, $K$-theory and derived equivalences, \emph{Duke
Math. J.} \textbf{124} (2004), no. 3, 587-617.}

\bibitem{GV}{{\sc S. Glaz} and {\sc W. V. Vasconcelos}, Flat ideals. II,\emph{ Manuscripta Math.}
\textbf{22} (1977), no. 4, 325-341.}

\bibitem{gkuk}{{\sc X. J. Guo} and {\sc A. Kuku}, Higher class groups of locally triangular
orders over number fields, \emph{Algebra Colloq.} \textbf{16}
(2009), no. 1, 79-84.}

\bibitem{hkx}{{\sc W. Hu, S. K\"onig, } and {\sc C. C. Xi},
Derived equivalences from cohomological approximations, and
mutations of $\Phi$-Yoneda algebras, Preprint,
arXiv:1102.2790v1[math.RT], 2011.}

\bibitem{hx2}{{\sc W. Hu } and {\sc C. C. Xi},
$\cal D$-split sequences and derived equivalences, \emph{Adv. Math.}
\textbf{227} (2011) 292-318.}

\bibitem{hx4}{{\sc W. Hu } and {\sc C. C. Xi}, Derived equivalences for $\Phi$-Auslander-Yoneda
algebras. Preprint, arXiv:9012.0647v2[math.RT], 2009.}

\bibitem{kaplansky}{{\sc I. Kaplansky}, \emph{Commutative rings}, University
of Chicago Press, Chicago, 1974.}

\bibitem{keating}{{\sc M. E. Keating}, The $K$-theory of triangular rings and orders, in:
\emph{Algebraic $K$-theory, number theory, geometry and analysis}
(Bielefeld, 1982), Lecture Notes in Math. 1046, 178-192,
Springer-Verlag, Berlin, 1984.}

\bibitem{kk}{{\sc E. Kirkman} and {\sc J. Kuzmanovich}, Global dimensions of a class of tiled orders, {\it J. Algebra} \textbf{127} (1989), no. 1, 57-72.}

\bibitem{kuk}{{\sc A. Kuku}, Higher algebraic $K$-theory for twisted Laurent series rings over
orders and semisimple algebras, \emph{Algebr. Represent Theor.}
\textbf{11} (2008) 355-368.}

\bibitem{magurn}{{\sc B. A. Magurn,} \emph{An algebraic introduction to $K$-theory}, Encyclopedia of Math. and its Applications 87, Cambridge University Press, 2002.}

\bibitem{milnor}{{\sc J. Milnor,} \emph{Introduction to algebraic $K$-theory}, Annals of Math. Studies 72, Princeton University Press, 1971.}

\bibitem{swan1}{{\sc R. G. Swan}, Algebric $K$-theory, \emph{Actes, Congr\`es Intern. Math.} \textbf{1} (1970) 191-199.}

\bibitem{Quillen}{{\sc D. Quillen}, Higher algebraic $K$-theory, I. In: \emph{Algebraic K-theory, I}: Higher $K$-theories
(Seattle, 1972), Lecture Notes in Math. 341, 85-147. Springer,
Berlin, 1973.}

\bibitem{qff}{{\sc D. Quillen}, On the cohomology and $K$-theory of the general linear groups over
a finite field, \emph{Ann. Math.} \textbf{96} (1972) 197-212.}

\bibitem{cr} {{\sc I. Reiner}, \emph{Maximal orders}, Academic Press, 1975.}

\bibitem{rickard}{{\sc J. Rickard}, Morita theory for derived categories,
\emph{J. London Math. Soc.} \textbf{39} (1989) 436-456.}

\bibitem{stafford}{{\sc J. T. Stafford,} Auslander-regular algebras and maximal orders, \emph{J. London Math. Soc.} \textbf{50} (1994) 276-292.}

\bibitem{wald}{{\sc F. Waldhausen}, Algebraic $K$-theory of generalized free products, I, II, \emph{Ann. Math.}\textbf{ 108}
(1978) 135-256.}

\bibitem{weibel1}{{\sc Ch. Weibel,} Mayer-Vietoris sequences and module structure on NK$_*$, In: \emph{Algebraic
K-theory} (Evanston, 1980), Lecture Notes in Math. 854, 466-493.
Springer-Verlag, Berlin, 1981.}

\bibitem{weibel2}{{\sc Ch. Weibel,} Mayer-Vietoris sequences and mod-p $K$-theory,
In: \emph{Algebraic $K$-theory}, Part I. Lecture Notes in Math. 966,
390-407. Springer-Verlag, Berlin, 1982.}

\bibitem{weibel}{{\sc Ch. Weibel,} \emph{An introduction to algebraic
$K$-theory}, A graduate textbook in progress.}
\bibitem{xi2}{{\sc C. C. Xi}, On the finitistic dimension conjecture, II: Related to finite global dimension, \emph{Adv. Math.} \textbf{201} (2006) 116-142.}

\bibitem{ywzc}{{\sc H. Y. Yin, F. G. Wang, X. S. Zhu} and {\sc Y. H. Chen}, $w$-modules over commutative rings, \emph{J. Korean Math. Soc.}\textbf{ 48}
(2011), no. 1, 207-222.}
\end{thebibliography}}
\end{document}